\documentclass[11pt]{article}
\usepackage{amssymb,amsmath,latexsym}
\usepackage{epsfig}

\newtheorem{thm}{Theorem}[section]

\newtheorem{cor}[thm]{Corollary}
\newtheorem{lemma}[thm]{Lemma}
\newtheorem{prop}[thm]{Proposition}

\newtheorem{remark}[thm]{Remark}
\numberwithin{equation}{section}

\def\1{{\bf 1}}
\def\ce{\mathcal{E}}
\def\cf{\mathcal{F}}

\def\pf{{\medskip\noindent {\bf Proof. }}}
\def\qed{{\hfill $\Box$ \bigskip}}

  \def\bR {{\mathbb R}}

\def\nn{\nonumber}

\def\wt{\widetilde}

\def\ol{\overline}
\def\E{{\mathbb E}}
\def\P{{\mathbb P}}

\def\bea{\begin{align*}}
\def\eea{\end{align*}}
\def\bee{\begin{equation}}
\def\eee{\end{equation}}

\def\vp{\varphi}

\def\R{{\mathbb R}}
\def\br{{\mathbb R}}
\def\b{{\beta}}

\def\ol{\overline}

\makeatletter
\@addtoreset{equation}{section}

\makeatother

\begin{document}
\allowdisplaybreaks
\bibliographystyle{plain}

\title{\Large \bf
Global heat kernel estimates for symmetric Markov processes dominated by stable-like processes in exterior $C^{1,\eta}$  open sets}

\author{{\bf Kyung-Youn Kim}\thanks{
This work was supported by the National Research Foundation of Korea(NRF) grant funded by the Korea government(MEST) (NRF- 2013R1A2A2A01004822).}
}

\date{}

\maketitle

\begin{abstract}
In this paper, we establish sharp two-sided heat kernel estimates for a large class of symmetric Markov processes in exterior  $C^{1,\eta}$  open sets for all $t> 0$. The processes are symmetric pure jump Markov processes with jumping kernel intensity
$$\kappa(x, y)\psi(|x-y|)^{-1}|x-y|^{-d-\alpha}$$
where $\alpha\in(0,2)$, $\psi$ is an increasing  function on $[ 0, \infty)$ with $\psi(r)=1$ on $0<r\le 1$ and $c_1e^{c_2r^{\b}}\le \psi(r)\le c_3e^{c_4r^{\b}}$ on $r>1$ for $\b\in[0, \infty]$. A symmetric function $\kappa(x, y)$ is bounded by two positive constants and $|\kappa(x, y)-\kappa(x,x)|\le c_5 |x-y|^{\rho}$ for $|x-y|<1$ and  $\rho>\alpha/2$. As a corollary of our main result, we estimates sharp two-sided Green function for this process in $C^{1,\eta}$ exterior open sets.
\end{abstract}


\bigskip

\section{Introduction}

In this paper, we study two-sided heat kernel estimates for a large class of symmetric Markov processes with jumps in   exterior $C^{1, \eta}$ open sets for all $t>0$. Discontinuous Markov processes and non-local Markovian operator have received much attention recently. The transition density $p(t, x, y)$ which describes the distribution of Discontinuous Markov process is a fundamental solution of involving infinitesimal generator and there are many studies in this areas in \cite{BBCK, CKK, CKK2,CKK3, CK, CK2}. Very recently in \cite{BGR2}, two-sided estimates on $p(t,x,y)$  for isotropic unimodal L\'evy processes with L\'evy exponents having weak local scaling at infinity are established. Also, heat kernel estimates for a class of L\'evy processes with L\'evy measures not necessarily absolutely continuous with respect to the underlying measure are obtained by Kaleta and Sztonyk in \cite{KaSz2}.

 Since it is difficult to obtain two-sided estimates on Dirichlet heat kernel where points are near the boundary, Dirichlet heat kernel estimates are obtained recently for particular processes in ~\cite{ BGR1, CKS, CKS2, CKS3}. Very recently, the studies of  two-sided Dirichlet heat kernel estimates are extended to a large class of symmetric L\'evy processes and beyond in~\cite{CKS8, CKS9, KK}.

In this paper, we consider a large class of symmetric Markov processes whose jumping kernels are dominated by the kernels of stable-like processes which  is discussed in \cite{KK}. Throughout this paper we assume that $\beta \in [0, \infty]$, $\alpha\in (0, 2)$, and $d \in\{1,2,3,\ldots\}$.
For two nonnegative functions $f$ and $g$, the notation $f\asymp g$ means that there are positive constants $c_1$ and 
$c_2$ such that $c_1g(x)\leq f (x)\leq c_2 g(x)$ in the common domain of definition for $f$ and $g$.
We will use the symbol ``$:=$,'' which is read as ``is defined to be.''

 Let $\psi$ be an increasing function on $[0, \infty )$ where $\psi(r)=1$ on $0 < r\leq 1$, and  $  L_1  e^{\gamma_1r^{\beta}} \leq \psi (r)\leq L_2  e^{\gamma_2r^{\beta}}$ on $ 1<r<\infty$. Here $L_1, L_2, \gamma_1,  \gamma_2$  are positive constants. For any $r>0$,  we define $ j(r) := r^{-d-\alpha}\psi(r)^{-1}$.
Let $\kappa(x,y)$ be a positive symmetric function which is satisfying
\begin{equation}\label{e:conkappa1}
 L_3^{-1}\le \kappa(x,y) \le L_3, \quad x,y \in \R^d,
\end{equation}
and for $\rho >\alpha/2$ ,
\begin{equation*}
 |\kappa(x,y)-\kappa(x,x)|{\bf 1}_{\{  |x-y|<1  \}}\leq L_4|x-y|^\rho, \quad x,y \in \R^d,
\end{equation*}
where  $L_3, L_4$ are positive constants.
We define a symmetric measurable function $J$ on $\br^d\times \br^d \setminus \{x=y\}$ as \begin{equation}\label{e:J2}
  J(x, y) :=  \kappa(x,y)   j(|x-y|) =
  \begin{cases}
  \kappa(x,y)  |x-y|^{-d-\alpha} \psi (|x-y|)^{-1}  & \text{ if  } \beta \in [0, \infty),\\
  \kappa(x,y)   |x-y|^{-d-\alpha}{\bf 1}_{\{|x-y| \le1\}} & \text{ if  } \beta = \infty.
  \end{cases}
\end{equation}

For any $u\in L^2(\br^d, dx)$, we define $\ce (u, u):=2^{-1}\int_{\br^d\times \br^d} (u(x)-u(y))^2 J(x, y) dx dy$ and $ \mathcal{D}(\ce):=\{f\in C_c (\br^d): \ce(f,f) <\infty\}$ where $C_c(\br^d)$ is the space of continuous functions with compact support in $\br^d$ equipped with uniform topology.
Let $\ce_1 (u, u):= \ce (u, u) +  \int_{\br^d} u(x)^2 dx$ and $\cf:=\overline{\mathcal{D}(\ce)}^{\ce_1}$. Then by \cite[Proposition 2.2]{CK2}, $(\ce, \cf)$  is a regular Dirichlet form on $L^2(\br^d,dx)$ and there is a Hunt process $Y$ associated with this on $\br^d$ (see~\cite{FOT}).

It is shown in \cite{KK} that the Hunt process $Y$ associated with $(\ce, \cf)$ is a subclass of the processes considered in~\cite{CKK3}.
Therefore, $Y$ is conservative and it has a H\"older continuous transition density $p(t, x, y)$ on $(0, \infty )\times \br^d\times \br^d$ with respect to the Lebesgue measure.  In \cite{KaSz1, Sz1}, this process is discussed and the upper bound estimates are obtained.

For any $x\in \bR^d$, stopping time $S$ with respect to the filtration of $Y$, and nonnegative measurable function $f$ on $\bR_+ \times \bR^d\times \bR^d$ where $f(s, y, y)=0$ for all $y\in\bR^d$ and $s\ge 0$, we have a  L\'evy system for $Y$ :
\begin{equation}\label{e:levy}
\E_x \left[\sum_{s\le S} f(s,Y_{s-}, Y_s) \right] = \E_x \left[ \int_0^S \left(\int_{\bR^d} f(s,Y_s, y) J(Y_s,y) dy \right) ds \right]
\end{equation}
(e.g., see \cite[Appendix A]{CK2}). It describes the jumps of the process $Y$, so the function $J$ is called the jumping intensity kernel of $Y$.

For $a, b\in \bR$, we use $\wedge$ and $\vee$ to denote $a\wedge b:=\min \{a, b\}$ and $a\vee b:=\max\{a, b\}$.
For any positive constants  $a, b, T$, we define functions $\Psi^1_{a, b, T}(t, r)$ on $(0, T] \times [0, \infty)$ as
\begin{align}\label{eq:qd}
\Psi^1_{a, b, T}(t, r):=
\begin{cases}
t^{-d /\alpha}\wedge  tr^{-d-\alpha}e^{-b r^\b}&\text{ if } \b\in[0,1],\\
t^{-d /\alpha}\wedge  tr^{-d-\alpha}&\text{ if } \b\in (1, \infty] \text{ with } r<1,\\
t \exp\left(-a \left( r \, \left(\log \frac{ T r}{t}\right)^{\frac{\b-1}{\b}}\wedge  r^\b\right) \right)\qquad &\text{ if } \b\in(1, \infty) \text{ with } r \ge 1,\\
\left(t/(T r)\right)^{ar}  &\text{ if }  \b=\infty \text{ with } r \ge 1
\end{cases}
\end{align}
and  $\Psi^2_{a, T}(t, r)$ on $[T, \infty)\times (0, \infty)$ as
\begin{align}\label{eq:qd1}
\Psi^2_{a, T}(t, r):=
\begin{cases}
t^{-d /\alpha}\wedge  tr^{-d-\alpha} &\text{if } \b=0,\\
t^{-d/2}\exp\left(-a \left(r^\b\wedge \frac{r^2}{t} \right)\right) &\text{if }\b\in(0, 1],\\
t^{-d/2}\exp\left(-a \left(r \left(1+\log^{+} \frac{Tr}{t} \right)^{(\b-1)/\b}\wedge\frac{ r^2}{t }\right)\right)& \text{if }\b\in(1, \infty), \\
t^{-d/2}\exp\left(-a \left(r\left(1+\log^{+}\frac{Tr}{t}\right)\wedge \frac{r^2}{t}\right)\right) &\text{if }\b=\infty
\end{cases}
\end{align}
where $\log^+x= \log x\cdot {\bf 1}_{ \{x\ge 1\}}+ 0\cdot {\bf 1}_{\{ x< 1\}}$.

By \cite[Theorem 1.2]{CK2}, \cite [Theorem 1.2 and Theorem 1.4]{CKK3} and \cite[Theorem 1.1]{KK}, it is known that for any $T>0$, there are positive constants $C_1, c\ge 1$ and $\gamma=\gamma(\gamma_1, \gamma_2)\ge 1$ such that
\begin{align}\label{eq:smT}
c^{-1}\Psi^1_{C_1, \gamma, T}(t, |x-y|)\,\le\,
p(t,x,y)\,\le\,
c\,\Psi^1_{C_1^{-1}, \gamma^{-1}, T}(t, |x-y|)
\end{align}
 for every $(t, x, y) \in (0,T]\times \R^d\times \R^d$ and
\begin{align}\label{eq:laT}
c^{-1}\Psi^2_{C_1, T}(t, |x-y|)\,\le\,
p(t,x,y)\,\le\,
c\,\Psi^2_{C_1^{-1}, T}(t, |x-y|)
\end{align}
for every $(t, x, y) \in [T, \infty)\times \R^d\times \R^d$. Even though in \cite[Theorem 1.2]{CK2} and \cite[Theorems 1.2 and 1.4]{CKK3} two-sided estimates for $p(t, x, y)$ are stated separately for the cases $0<t\le 1$ and $t\ge1$, the constant 1 does not play any special role.
Thus by the same proof, two-sided estimates for $p(t, x, y)$ hold for the case $0<t\le T$ and can be stated in the above way. We remark here that in \cite[Theorems 1.2(2.b)]{CKK3} the case $|x-y|\asymp t$ is missing.
One can see that  \eqref{eq:laT} is the correct form to include the case $|x-y|\asymp t$ (cf. Proposition \ref{step4} below for the lower bound).

The goal of this paper is to establish the two-sided heat kernel estimates for $Y$ in exterior $C^{1,\eta}$ open set. Recall that an open set $D$ in $\bR^d$ (when $d\ge 2$) is said to be $C^{1,\eta}$ open set   with $\eta\in(0, 1]$ if there exist a localization radius $r_0>0 $ and a constant $\Lambda_0>0$ such that for every $z\in\partial D$, there exists a $C^{1,\eta}$-function $\phi=\phi_z: \bR^{d-1}\to \bR$ satisfying $\phi(0)=0$, $\nabla\phi (0)=(0, \dots, 0)$, $\| \nabla\phi \|_\infty \leq \Lambda_0$, $| \nabla \phi (x)-\nabla \phi (w)| \leq \Lambda_0 |x-w|^{\eta}$ and an orthonormal coordinate system $CS_z$ of  $z=(z_1, \cdots, z_{d-1}, z_d)=:(\wt z, \, z_d)$ with origin at $z$ such that $B(z, r_0 )\cap D= \{y=({\tilde y}, y_d) \in B(z, r_0) \mbox{ in } CS_z: y_d > \phi (\wt y) \}$.
The pair $(r_0, \Lambda_0)$ will be called the $C^{1,\eta}$ characteristics of the open set $D$.
Note that a $C^{1,\eta}$ open set $D$ with characteristics $(r_0, \Lambda_0)$ can be unbounded and disconnected.

Let $Y^D$ be the subprocess of $Y$ killed upon exiting $D$ and $\tau_D:=\inf\{t>0: Y_t\notin D\}$ be the first exit time from $D$.  By the strong Markov property, it can easily be verified that
$p_D(t, x, y):=p(t, x, y)-\E_x[ p(t-\tau_D, Y_{\tau_D}, y); t>\tau_D]$ is the transition density of $Y^D$. Also, by the continuity and estimate of $p$, it is routine to show that $p_D(t, x, y)$ is symmetric and continuous(e.g., see the proof of  Theorem 2.4 in \cite{CZ}).

In \cite [Theorem 1.2]{KK}, the Dirichlet heat kernel estimates for $Y^D$ is obtained. For the lower bound estimates on $p_D(t,x,y)$ when $\beta \in (1, \infty]$, we need the following assumption on $D$: {\it the path distance in each connected component of $D$ is comparable to the Euclidean distance with characteristic $\lambda_1$}, i.e., for every $x$ and $y$ in the same component of $D$ there is a rectifiable curve $l$ in $D$ which connects $x$ to $y$ such that the length of $l$ is less than or equal to  $\lambda_1|x-y|$.
Clearly, such a property holds for all bounded $C^{1,\eta}$ open sets, $C^{1,\eta}$ open sets with compact complements, and connected open sets above graphs of $C^{1,\eta}$ functions.

Here is the main result of \cite{KK}. We denote by $\delta_D(x)$ the Euclidean distance between $x$ and $D^c$.

\begin{thm}\label{T:sm} \emph{\cite [Theorem 1.2]{KK}}
Let $J$ be the symmetric function defined in~\eqref{e:J2} and $Y$ be the symmetric pure jump Hunt process with the jumping intensity kernel $J$. Suppose that $T>0$ and $\gamma$ is the constant in  \eqref{eq:smT}.
For any $\eta \in (\alpha/2,1]$,  let $D$ be a $C^{1,\eta}$ open set in $\R^d$ with $C^{1, \eta}$ characteristics $(r_0, \Lambda_0)$. Then the transition density $p_D(t, x, y)$ of $Y^D$ has the following estimates.
\begin{description}
\item{\rm (1)}
There are positive constants $c, C_2\ge 1$ such that for any $(t, x, y)\in(0, T]\times D\times D$, we have
\begin{align*}
& c\, \left(1\wedge\frac{\delta_D(x)}{t^{1/\alpha}}\right)^{\alpha/2}\,\left(1\wedge\frac{\delta_D(y)}{t^{1/\alpha}}\right)^{\alpha/2} \Psi^1_{C_2^{-1}, \gamma^{-1},  T}(t, |x-y|/6) \ge p_D(t, x, y)\nn\\
&\ge c^{-1} \, \left(1\wedge\frac{\delta_D(x)}{t^{1/\alpha}}\right)^{\alpha/2}\left(1\wedge\frac{\delta_D(y)}{t^{1/\alpha}}\right)^{\alpha/2} \cdot\begin{cases}
t^{-d/\alpha}\wedge t|x-y|^{-d-\alpha} e^{-\gamma |x-y|^\b}&\hbox{if } \b\in[0,1],\\
t^{-d/\alpha}\wedge t|x-y|^{-d-\alpha}&
\begin{array}{c}
\hbox{ \hskip -.1in  if  $\b\in(1, \infty]$ and }\\
\hbox{\hskip -.3in $|x-y|\le 4/5$}.
\end{array}
\end{cases}
\end{align*}
\item{\rm(2)}
Suppose in addition that the path distance in each connected component of $D$ is comparable to the Euclidean distance with characteristic $\lambda_1$. If $\b\in (1,\infty]$,  there are positive constants $c, C_2\ge 1$ such that for any $(t, x, y)\in (0, T]\times D\times D$ where $|x-y|\ge 4/5$ and $x, y$ are in a same component of $D$, we have
\begin{align*}
p_D(t, x, y) \ge  c^{-1} &  \left(1\wedge\frac{\delta_D(x)}{t^{1/\alpha}}\right)^{\alpha/2}\,\left(1\wedge\frac{\delta_D(y)}{t^{1/\alpha}}\right)^{\alpha/2}\Psi^1_{C_2,  \gamma, T}(t, 5|x-y|/4)
\end{align*}

\item{\rm (3)}
If $\b\in (1, \infty)$, there is a positive constant $c\ge 1$ such that for any $(t, x, y)\in (0, T]\times D\times D$ where $|x-y|\ge 4/5$ and $x, y$ are in different components of $D$, we have
\begin{align*}
  p_D(t, x, y) \ge\,
c^{-1}\,    \left(1\wedge\frac{\delta_D(x)}{t^{1/\alpha}}\right)^{\alpha/2}&\,\left(1\wedge\frac{\delta_D(y)}{t^{1/\alpha}}\right)^{\alpha/2}\frac{t}{ |x-y|^{d+\alpha}}e^{-\gamma (5|x-y|/4)^\b}.
\end{align*}

\item{\rm (4)}
Suppose in addition that $D$ is bounded and connected.
Then there is positive constant $c\ge 1$ such that for any  $(t, x, y)\in [T, \infty)\times D\times D$ we have
\[
c^{-1}\, e^{- t \, \lambda^{D} }\, \delta_D (x)^{\alpha/2}\, \delta_D (y)^{\alpha/2}\,\leq\,
p_D(t, x, y) \,\leq\,
c\,e^{- t\, \lambda^{D}}\, \delta_D (x)^{\alpha/2} \,\delta_D (y)^{\alpha/2},
\]
where $-\lambda^D<0$ is the largest eigenvalue of the generator of $Y^D$.
\end{description}
\end{thm}

Theorem \ref{T:sm}(1)--(3) give us the Dirichlet heat kernel estimates for the small time. However the large time estimates are established only for the bounded and connected $C^{1, \eta}$ open sets. The large time Dirichlet heat kernel estimates for unbounded open sets are different depending on the geometry of $D$ as one sees for the cases of the symmetric $\alpha$-stable processes and of the relativistic stable processes in \cite{CT} and in \cite{ CKS6, CKS7}, respectively.

Motivated by \cite{CT, CKS7}, we establish the global sharp two-sided estimates on $p_D(t, x, y)$ in the exterior $C^{1,\eta}$ open set, that is, $C^{1,\eta}$ open set which is $D^c $ is compact.
It can be disconnected and in this case, there are bounded connected components.
The number of the such bounded connected components is finite.

\begin{thm}\label{T:main}
Let $J$ be the symmetric function defined in~\eqref{e:J2} and $Y$ be the symmetric pure jump Hunt process with the jumping intensity kernel $J$.  Let  $d>2\cdot{\bf 1}_{\{ \beta \in (0, \infty]\}}   + \alpha\cdot {\bf 1}_{\{ \beta =0\}} $, $T>0$ and $R>0$ be positive constants.
For any $\eta \in (\alpha/2,1]$,  let $D$  be an exterior $C^{1, \eta}$ open set in $\R^d$ with $C^{1, \eta}$ characteristics $(r_0, \Lambda_0)$ and $D^{\,c}\subset B(0, R)$.
Let $D_0$ be an unbounded connected component and $D_1,\ldots, D_n$ be bounded connected components such that $D_0\cup D_1\cup\ldots \cup D_n=D$.
Then for any $t\ge T$ ($t>0$ when $\b=0$, respectively) and $ x, y\in D$,
the transition density $p_D(t, x, y)$ of $Y^D$ has the following estimates.
\begin{description}
\item{\rm (1)}
For any $\b\in[0, \infty]$,  there are positive constants $c_i=c_i(\alpha, \beta, \eta, r_0, \Lambda_0, R, T, d, L_3, L_4, \psi)$ ($c_i=c_i(\alpha, \beta, \eta, r_0, \Lambda_0, R, d, L_3, L_4, \psi)$ when $\b=0$, respectively), $i=1,2$  such that
\begin{align*}
p_D(t, x, y)\le c_1 &\left(1\wedge \frac{ \delta_D(x)}{1\wedge t^{1/\alpha}}\right)^{\alpha/2}\left(1\wedge \frac{ \delta_D(y)}{1\wedge t^{1/\alpha}} \right)^{\alpha/2}\Psi^2_{c_2, T}(t, |x-y|).
\end{align*}
\item{\rm (2)}
Suppose that $\b\in[0, 1]$ or  $\b\in(1, \infty]$ with $|x-y|< 4/5$. Then  there are positive constants $c_i=c_i(\alpha, \beta,  \eta, r_0, \Lambda_0, R, T, d,L_3, L_4, \psi)$ ($c_i=c_i(\alpha, \beta, \eta, r_0, \Lambda_0, R, d, L_3, L_4, \psi)$ when $\b=0$, respectively), $i=1,2$  such that
\begin{align*}
p_D(t, x, y) &\ge c_1 \left(1\wedge \frac{ \delta_D(x)}{1\wedge t^{1/\alpha}}\right)^{\alpha/2}\left(1\wedge \frac{ \delta_D(y)}{1\wedge t^{1/\alpha}} \right)^{\alpha/2}\Psi^2_{c_2, T}(t, |x-y|).
\end{align*}
\item{\rm (3)}
Suppose that $\b\in(1, \infty]$ with $|x-y|\ge 4/5$ and $x, y$ are in a same component of $D$.
\begin{description}
\item{\rm(3.a)}  \emph{(Unbounded connected component)}
There are positive constants $c_i=c_i(\alpha, \beta, \eta, r_0, \Lambda_0, $ $R, T, d,$ $ L_3, L_4, \psi)$, $i=1,2$  such that for $x,y \in D_0$
\begin{align*}
p_D(t, x, y)\,\ge \,c_1 \left(1\wedge \delta_D(x)\right)^{\alpha/2}\left(1\wedge \delta_D(y)\right)^{\alpha/2}
\Psi^2_{c_2, T}(t, |x-y|).
\end{align*}
\item{\rm(3.b)}
 \emph{(Bounded connected component)}
There is a  positive constant $c=c(\alpha, \beta, \eta, r_0, \Lambda_0, $ $R, T, d, L_3, L_4, \psi)$ such that   if  $x,y \in D_j$ for some $j=1, \dots n$,
 \begin{align*}
p_D(t, x, y) \,\ge\,
c\,e^{- t\, \lambda_j}\, \delta_D (x)^{\alpha/2} \,\delta_D (y)^{\alpha/2}
\end{align*}
where $-\lambda_{j}<0$ is the largest eigenvalue of the generator $Y^{D_j}$, $j=1\ldots, n$.
\end{description}
\item{\rm (4)}
Suppose that  $\b\in(1, \infty)$ with $|x-y|\ge 4/5$ and $x, y$ are in different components of $D$.
Then  there are positive constants $c_i=c_i(\alpha, \beta, \eta, r_0, \Lambda_0, R, T, d,L_3, L_4, \psi, \lambda_1,\ldots,\lambda_n)$, $i=1,2$ such that
\begin{align*}
p_D(t, x, y)\ge c_1 \left(1\wedge \delta_D(x)\right)^{\alpha/2}\left(1\wedge \delta_D(y)\right)^{\alpha/2}\frac{\exp\left(-c_2(|x-y|^{\b}+t)\right)}{|x-y|^{d+\alpha}}
\end{align*}
where $-\lambda_{j}<0$ is the largest eigenvalue of the generator $Y^{D_j}$, $j=1\ldots, n$.
\end{description}
\end{thm}

For a connected exterior $C^{1, \eta}$ open set, we can rewrite the sharp two-sided estimates on $p_D(t, x, y)$  for all $t> 0$ in a simple form combining Theorem \ref{T:sm}(1)--(2) and Theorem \ref{T:main}(1)--(3a).

\begin{cor}\label{Coro}
Let $J$ be the symmetric function defined in~\eqref{e:J2} and $Y$ be the symmetric pure jump Hunt process with the jumping intensity kernel $J$. Let  $d>2\cdot{\bf 1}_{\{ \beta \in (0, \infty]\}}   + \alpha\cdot {\bf 1}_{\{ \beta =0\}} $, $T>0$ and $R>0$ be positive constants.
For any $\eta \in (\alpha/2,1]$,  let $D$  be a connected exterior $C^{1, \eta}$ open set in $\R^d$ with $C^{1, \eta}$ characteristics $(r_0, \Lambda_0)$ and $D^{\,c}\subset B(0, R)$.
Then there are positive constants $c_i=c_i(\alpha, \beta, \eta, r_0, \Lambda_0, R, T, d, L_3, L_4, \psi)>1, i=1,2$  such that for every $(t, x, y)\in(0, \infty)\times D\times D$, we have
\begin{align*}
p_D(t, x, y)\le c_1 \left(1\wedge \frac{\delta_D(x)}{1\wedge t^{1/\alpha}}\right)^{\alpha/2}\left(1\wedge \frac{\delta_D(y)}{1\wedge t^{1/\alpha}}\right)^{\alpha/2}
\cdot
\begin{cases}
\Psi^1_{c_2^{-1},\gamma^{-1}, T}(t, |x-y|/6)&\hbox{if } t\in(0,T],\\
\Psi^2_{c_2^{-1}, T}(t, |x-y|)&\hbox{if } t\in[T,\infty),
\end{cases}
\end{align*}
and in addition $D$  is a connected, we have
\begin{align*}
p_D(t, x, y)\ge  c_1^{-1} \left(1\wedge \frac{\delta_D(x)}{1\wedge t^{1/\alpha}}\right)^{\alpha/2}\left(1\wedge \frac{\delta_D(y)}{1\wedge t^{1/\alpha}}\right)^{\alpha/2}
\cdot
\begin{cases}
\Psi^1_{c_2,\gamma, T}(t, 5|x-y|/4)&\hbox{if } t\in(0,T],\\
\Psi^2_{c_2, T}(t, |x-y|)&\hbox{if } t\in[T,\infty)
\end{cases}
\end{align*}
where $\gamma$ is the constant in Theorem \ref{T:sm}.
\end{cor}

By integrating the heat kernel estimates in  Corollary \ref{Coro} with respect to $t\in(0, \infty)$, one gets the following sharp two-sided Green function estimates of $Y^D$ in the connected exterior $C^{1,\eta}$ open sets.

 \begin{cor}\label{C:green}
Let $J$ be the symmetric function defined in~\eqref{e:J2} and $Y$ be the symmetric pure jump Hunt process with the jumping intensity kernel $J$.   Let  $d>2\cdot{\bf 1}_{\{ \beta \in (0, \infty]\}}   + \alpha\cdot {\bf 1}_{\{ \beta =0\}} $ and $R>0$ be a positive constant.
For any $\eta \in (\alpha/2,1]$,  let $D$  be a connected exterior $C^{1, \eta}$ open set in $\R^d$ with $C^{1, \eta}$ characteristics $(r_0, \Lambda_0)$ and $D^{\,c}\subset B(0, R)$.
Then there is a positive constant $c=c(\alpha,\b, \eta,  r_0, \Lambda_0, R, d , L_3, L_4,\psi)>1$ such that  for every $(x, y)\in D\times D$, we have
\begin{align*}
&c^{-1} \left(\frac{1}{|x-y|^{d-\alpha}}+\frac{1}{|x-y|^{d-2}}\cdot{\bf 1}_{\{ \beta \in (0, \infty]\}}\right)\left(1\wedge \frac{\delta_D(x)}{|x-y|\wedge 1}\right)^{\alpha/2}\left(1\wedge \frac{\delta_D(y)}{|x-y|\wedge 1}\right)^{\alpha/2}\\
\le&~ G_D(x, y)\le c\left(\frac{1}{|x-y|^{d-\alpha}}+\frac{1}{|x-y|^{d-2}}\cdot{\bf 1}_{\{ \beta \in (0, \infty]\}}\right)\left(1\wedge \frac{\delta_D(x)}{|x-y|\wedge 1}\right)^{\alpha/2}\left(1\wedge \frac{\delta_D(y)}{|x-y|\wedge 1}\right)^{\alpha/2}.
\end{align*}
\end{cor}

The  approach developed in \cite{CKS7} provides us
a main road map.
By checking the cases depending on  the value of $\beta$ and  the distance between $x$ and $y$ carefully,  we establish sharp two-sided estimates on $p_D(t, x, y)$ for exterior $C^{1,\eta}$ open sets for all $t\in [T, \infty)$.
In section 2, we first give elementary results on the functions $\Psi^1(t, r)$ and $\Psi^2(t, r)$  which are defined in \eqref{eq:qd} and \eqref{eq:qd1}. Also, we give the proof of the upper bound estimates on $p_D(t, x, y)$. In Section 3, we present the interior lower bound estimates on $p_{\ol{B}_R^c}(t, x, y)$ where $B_R:=B(x_0, R)$ for some $x_0\in \R^d$. In Section 4, the full lower bound estimates on $p_D(t, x, y)$  for exterior  open set $D$ are established by considering the cases whether the points are in a same component or in different components separately. The proof of Corollary \ref{C:green} is given in Section 5.

Throughout this paper, the positive constants $C_1, C_2,$ $L_1, L_2, L_3, L_4, \gamma_1, \gamma_2, \gamma$ will be fixed.  In the statements of results and the proofs, the constants $c_i=c_i(a,b,c,\ldots), i=1,2,3,\ldots$, denote generic constants depending on $a, b, c, \ldots$ and there are given anew in each statement and each proof. The dependence of the constants on the dimension $d$, on $\alpha\in(0,2)$ and on the positive constants $L_1, L_2, L_3, L_4,\gamma_1, \gamma_2, \gamma$ will not be mentioned explicitly.

\section{Upper bound estimates}

We first give elementary lemmas which are used several times to estimates the upper and lower bound on $p_D(t, x, y)$ where $t\ge T$ ($t>0$ when $\b=0$, respectively). Recall the functions $\Psi^1(t, r)$ and $\Psi^2(t, r)$ which  are  defined in \eqref{eq:qd} and \eqref{eq:qd1}.

\begin{lemma}\label{L:p1}
Let $t_0>0$ and $a, b, c\ge 1$ be fixed constants. For any $\b\in (0, \infty]$, suppose that  $N_1, N_2$ be positive constants satisfying $N_2\ge N_1 \cdot (ab\vee c^{2/\b})$.
Then there exist positive constants $c_i=c_i(t_0), i=1,2$ such that  for every  $r>0$, we have that
\begin{align*}
(1)&\qquad\Psi^1_{b^{-1}, c^{-1}, t_0}(t_0, N_1^{-1}r)\le c_1\Psi^1_{a, c, t_0}(t_0, N_2^{-1}r)\qquad\mbox{and}\nn\\
(2)&\qquad\Psi^1_{a^{-1}, c^{-1}, t_0}(t_0, N_2r)\le c_2\Psi^1_{b, c, t_0}(t_0, N_1r).
\end{align*}
\end{lemma}
\pf
When $\b\in(0,1]$, since $N_2\ge N_1 c^{2/\b}$, we have (1) and (2).

When $\b\in(1,\infty]$,  since $t_0^{-d/\alpha}\wedge t_0r^{-d/\alpha}\asymp 1$ for any $r<1$, we only consider the case $1\le N_2^{-1}r(\le  N_1^{-1}r)$ to prove (1) and  $1\le N_1r(\le N_2r)$ to prove (2).
In these cases, since $\log x$ is increasing in $x$ and  $N_2\ge N_1ab$, we have (1) and (2).
\qed

\begin{lemma}\label{L:p2}
Let $T,a$ and $ b$ be positive constants. (1) If $b\ge 1$, there exists a positive constant $c=c(b)$ such that  for every $t\in[T, \infty)$ and $ r>0$, we have that
\begin{align*}
&\Psi^2_{a, T}(t, b^{-1}r)\le \Psi^2_{ab^{-2}, T}(t, r).
\end{align*}
(2) In addition,  for $a, b \ge 1$ and $\b\in (0, \infty]$, suppose that $N$ be a positive constant satisfying $N\ge (ab)^{1/(\b\wedge 1)}$. Then for every $t\in[T, \infty)$ and $ r>0$, we have that
\begin{align*}
&\Psi^2_{b^{-1}, T}(t, r)\le \Psi^2_{a, T}(t, N^{-1}r).
\end{align*}
\end{lemma}
\pf
Since $b\ge 1$, it is easy to prove (1) when  $\b\in[0, 1]$.
Also, since
$$b\left(1+\log^+\frac{Tb^{-1}r}{t}\right)\ge\left(1+\log b\right)\cdot\left(1+\log^+\frac{Tb^{-1}r}{t}\right)\ge \left(1+\log^+\frac{Tr}{t}\right),$$
 for any $b\ge1$, we have (1) when $\b\in(1, \infty]$.

On the other hand,  since $N\ge (ab)^{1/\b}(\ge 1)$, we have that
\begin{align}\label{eq:b_1}
b^{-1}\left(r^\b\wedge \frac{r^2}{t}\right)&\ge b^{-1}N^{\b}\left((N^{-1}r)^\b\wedge \frac{(N^{-1}r)^2}{t}\right)\ge a\left((N^{-1}r)^\b\wedge \frac{(N^{-1}r)^2}{t}\right).
\end{align}
Also, since $r\to 1+\log^+r$ is non-decreasing and $N\ge ab\,(\ge 1)$, we have that
\begin{align}\label{eq:b_2}
b^{-1}\left(r\left(1+\log^+\frac{Tr}{t}\right)^{(\b-1)/\b}\wedge \frac{r^2}{t}\right)&\ge b^{-1}N\left(N^{-1}r\left(1+\log^+\frac{N^{-1}Tr}{t}\right)^{(\b-1)/\b}\wedge \frac{(N^{-1}r)^2}{t}\right)\nn\\
&\ge a\left(N^{-1}r\left(1+\log^+\frac{N^{-1}Tr}{t}\right)^{(\b-1)/\b}\wedge \frac{(N^{-1}r)^2}{t}\right).
\end{align}
Hence, by \eqref{eq:b_1} for $\b\in(0, 1]$ and by \eqref{eq:b_2} for $\b\in(1, \infty]$, we have (2).
\qed

We now prove the upper bound estimates in Theorem \ref{T:main}(1).

\medskip
\noindent {\bf Proof of  Theorem \ref{T:main}(1)}
When $\b=0$, by Theorem \ref{T:sm}(1), we may assume that $t\ge T$.
Wthout loss of the generality, we may assume that $T=3$.
By the semigroup property and Theorem \ref{T:sm}(1), we have that for $t-2\ge 1$ and $x, y\in D$,
\begin{align}\label{e:up}
p_D(t, x, y)=&\int_D\int_Dp_D(1, x, z)p_D(t-2, z, w)p_D(1, w, y)dzdw\nn\\
\le&\, c_1 \left(1\wedge \delta_D(x)\right)^{\alpha/2}\left(1\wedge \delta_D(y)\right)^{\alpha/2} f_1(t, x, y)
\end{align}
where  $C_2$ and $\gamma$ are given constants in Theorem \ref{T:sm} and
\begin{align}\label{e:f_1}
f_1(t, x, y)=\int_{\R^d\times \R^d}\Psi^1_{C_2^{-1}, \gamma^{-1}, 1}(1, |x-z|/6) \,p(t-2, z, w)\, \Psi^1_{C_2^{-1}, \gamma^{-1}, 1}(1, |y-w|/6)dzdw.
\end{align}

Let $A_1:=\max\{ C_1^{2/(\b\wedge 1)}, 6\gamma^{2/\b},  6C_1C_2\}$ ($A_1= 6$ when $\b=0$, respectively) where $C_1$ is given constant in \eqref{eq:smT} and  \eqref{eq:laT}.
Then by \eqref{eq:laT},  there exists constants $c_i=c_i(\b)>0$, $i=2,3$ such that
\begin{align*}
p(t-2, z, w)\le c_2\,\Psi^2_{C_1^{-1}, 1}(t-2, |z-w|)\le c_2\,\Psi^2_{C_1, 1}(t-2, A_1^{-1}|z-w|)\le c_3\,p(t-2,  A_1^{-1}z, A_1^{-1}w).
\end{align*}
For the second inequaltiy, when $\b\in(0, \infty]$, we use (2) in Lemma \ref{L:p2} with $N=A_1$,  $a=b=C_1$ and the fact $A_1\ge C_1^{2/(\b\wedge 1)}$. When $\b=0$, the second inequality holds since $A_1\ge 1$.

Also, by \eqref{eq:smT},  there exist constants $c_i=c_i(\b)>0, i=4,5$ such that
\begin{align*}
\Psi^1_{C_2^{-1}, \gamma^{-1}, 1}(1, |x-z|/6)&\le c_4\, \Psi^1_{C_1 \gamma, 1}(1, A_1^{-1}|x-z|)\le c_5\, p(1, A_1^{-1}x, A_1^{-1}z)\,\,\mbox{ and } \\
\Psi^1_{C_2^{-1}, \gamma^{-1}, 1}(1, |y-w|/6)&\le c_4 \, \Psi^1_{C_1 \gamma, 1}(1, A_1^{-1}|y-w|)\le c_5\, p(1, A_1^{-1}y, A_1^{-1}w).
\end{align*}
For the first inequalties above, when $\b\in(0, \infty]$, we use
(1) in Lemma \ref{L:p1} along with $a=C_1$, $b=C_2$, $c=\gamma$, $N_1=6$ and  $N_2=A_1$ and the fact $A_1\ge 6( C_1C_2\vee \gamma^{2/\b})$. When $\b=0$, the first inequalities hold since $A_1= 6$.

Applying the above observations to \eqref{e:f_1} and  by the change of variable $\hat{z}=A_1^{-1}z$,  $\hat{w}=A_1^{-1}w$,  the semigroup property and \eqref{eq:laT}, we conclude that
\begin{align}\label{f_1}
f_1(t, x, y)\le& ~c_6 \int_{\R^d\times \R^d} p(1, A_1^{-1}x, \hat{z})\, p(t-2, \hat{z}, \hat{w} ) p(1, A_1^{-1}y, \hat w)d\hat{z} d\hat{w}\nn\\
=&~ c_6\, p(t, A_1^{-1}x, A_1^{-1}y)\le c_7\Psi^2_{C_1^{-1}, T}(t, A_1^{-1}|x-y|)\nn\\
\le&~ c_8\Psi^2_{C_1^{-1}A_1^{-2}, T}(t, |x-y|).
\end{align}
We have applied (1) in  Lemma \ref{L:p2} with $a=C_1$ and $b=A_1$  for the last inequality. Applying \eqref{f_1} to  \eqref{e:up}, we have proved the upper bound estimates in Theorem \ref{T:main}.
\qed

\section{Interior lower bound estimates}
The goal of this section is to the establish interior lower bound estimate on the heat kernel $p_{\ol{B}_R^c}(t, x, y)$ for  $t\ge T$  ($t>0$ when $\b=0$, respectively) where $B_R=B(x_0, R)$ for some $R>0$ and $x_0\in \R^d$.  We will combine ideas from \cite{CKS6} and \cite{KK}.

First, we introduce a Lemma which will be used in the proof of Lemma \ref{L:3.3} and Proposition \ref{step1}. Let $\vp(r):=r^2\cdot {\bf 1}_{\{\beta\in (0, \infty]\}}+r^{\alpha}\cdot{\bf 1}_{\{\beta=0\}}$ and then $\vp^{-1}(t)=t^{1/2}\cdot {\bf 1}_{\{\beta\in (0, \infty]\}}+t^{1/\alpha}\cdot{\bf 1}_{\{\beta=0\}}$.

\begin{lemma}\label{L:3.2}
 Let $a$ be a positive constant and  $T>0$ and $\b\in[0,\infty]$. Then there exists a positive constant $c=c(a, \b, T)$   ($c=c(a)$ when $\beta=0$, respectively) such that  for all $t \in [T, \infty)$ ($t>0$ when $\beta=0$, respectively), we have
\[
\inf_{y\in\bR^d  } \P_y \left(\tau_{B(y, a\vp^{-1}(t))} > t \right)\, \ge\, c_.
\]
\end{lemma}

\pf
When $\b=0$, using \cite[Theorem 4.12 and Proposition 4.9]{CK2},  the proof is almost identical to  that of \cite  [Lemma 3.1] {CKS3}. When $\beta\in(0, \infty]$, using \cite[Theorem 4.8]{CKK3}, the proof is the same as that of \cite  [Lemma 3.2] {CKS6}. So we omit the proof detail.
\qed

\begin{lemma}\label{L:3.3}
Let $D$ be an arbitrary open set. Suppose that $a$  be a positive constant and $T>0$ and $\b\in[0,\infty]$.
Then there exists a positive constant $c=c(a,\b, T)$   ($c=c(a)$ when $\b=0$, respectively) such that  for all $t\in [T,\infty)$  ($t>0$ when $\beta=0$, respectively) and $x, y\in D$  with $\delta_D(x)\wedge \delta_D (y) \ge a \vp^{-1}(t)$ and $|x-y|\ge 2^{-1} a \vp^{-1}(t)$,  we have
\[
\P_x \left( Y^D_t \in B \big( y, \,  2^{-1} a \vp^{-1}(t) \big) \right)\,\ge c\,  t\cdot \vp^{-d}(t)j(|x-y|).
\]
\end{lemma}

\pf
Using Lemma \ref{L:3.2}, the strong Markov property and  L\'evy system \eqref{e:levy}, the proof of the lemma is similar to that of   \cite[Proposition 3.3]{KK}. So we omit the proof detail.
\qed

For the remainder of this section,  we assume that $D$ is a domain with the following property: there exist $\lambda_1 \in [1, \infty)$ and $\lambda_2 \in (0, 1]$ such that for every $r \le 1$ and  $x,y$ in the same component of $D$ with $\delta_D(x)\wedge \delta_D(y)\ge r$, there exists in $D$ a length parameterized rectifiable curve $l$ connecting $x$ to $y$ with the length $|l|$ of $l$ is less than or equal to $\lambda_1|x-y|$ and $\delta_D(l(u))\geq\lambda_2 r$ for $u\in(0,|l|].$
Clearly, such a property holds for all  $C^{1, \eta}$ domains with compact complements, and domains above graphs of $C^{1,\eta}$ functions.

The following Propositions are motivated by  \cite{CKS6}.

\begin{prop}\label{step1}
Let $a$  be a positive constant and $T>0$ and $\b\in[0,\infty]$.
Then there exists a positive constant $c=c(a,\b, T, \lambda_1, \lambda_2)$ ($c=c(a, \lambda_1, \lambda_2)$ when $\b=0$, respectively) such that  for all $t\in [T,\infty)$ ($t>0$ when $\beta=0$, respectively) and $x, y\in D$  with $\delta_D(x)\wedge \delta_D (y) \ge a \vp^{-1}(t)\ge 2|x-y|$,  we have $ p_D(t,x,y) \,\ge\,c\, /\vp^{-d}(t).$
\end{prop}
\pf
By the same proof as that of \cite[Proposition 3.4]{CKS6}, we deduce the proposition using the parabolic Harnack inequality(see \cite[Theorem 4.12]{CK2} for $\b=0$ and  \cite [Theorem 4.11]{CKK3} for $\b\in(0, \infty]$) and Lemma \ref{L:3.3}.
\qed

\begin{prop}\label{step2}
Let  $a$  be a positive constant and $T>0$ and $\b\in[0,\infty]$.
Then there exists a positive constant $c=c(a,\b, T, \lambda_1, \lambda_2)$ ($c=c(a, \lambda_1, \lambda_2)$ when $\b=0$, respectively) such that  for all $t\in [T,\infty)$ ( $t>0$ when $\beta=0$, respectively) and $x, y\in D$  with $\delta_D(x)\wedge \delta_D (y) \ge a \vp^{-1}(t)$ and $|x-y|\ge 2^{-1} a \vp^{-1}(t)$,  we have $p_D(t, x, y)\,\ge c\, t j(|x-y|)$.
\end{prop}
\pf
By the same proof as that of \cite[Proposition 3.5]{CKS6}, we deduce the proposition using  the semigroup property, Lemma \ref{L:3.3} and Proposition \ref{step1}.
\qed

Also, since the  proof of the following proposition is almost identical to that of  \cite[Proposition 3.6]{CKS6} using Proposition \ref{step1}, we skip the proof.

\begin{prop}\label{step3}
Let $\b\in(1, \infty]$ and $a$ and $C_*$ be positive constants. Then there exist positive constants $c_i=c_i(a, \b,  C_*,\lambda_1, \lambda_2)$, $i=1,2$  such that  for every $t\in (0, \infty)$ and $x, y\in D$ with $\delta_D(x)\wedge \delta_D(y) \ge a\sqrt{ t}$, we have
\begin{align*}
p_D(t,x,y) \,\ge\,c_1\, t^{-d/2}\exp\left(-c_2\frac{|x-y|^2}{t}\right) ~\mbox{when }\,  C_*|x-y|\le t \le |x-y|^2.
\end{align*}
\end{prop}

Now, we estimates the interior lower bound for $p_D(t, x, y)$ where $\b\in(1, \infty]$ and $T\le t\le C_*T|x-y|$ for any positive constant $C_*<1$. The following Proposition \ref{step4} and Proposition \ref{step5} are counterparts of \cite[Propsition 3.6]{KK} and \cite[Proposition 3.5]{KK}, respectively. (See, also \cite[Theorem 5.5]{CKK3}) and \cite[Theorem 3.6]{CKK},  respectively.)

\begin{prop}\label{step4}
Let $\b\in(1, \infty)$ and $a, T$ and $C_*\in(0,1)$ be positive constants. Then there exist positive constants $c_i=c_i(a, \b,T,  C_*,\lambda_1, \lambda_2)$, $i=1,2$  such that  for every $t\in [T, \infty)$ and $x, y\in D$ with $\delta_D(x)\wedge \delta_D(y) \ge a\sqrt{ t}$, we have
\[
p_D(t,x,y) \,\ge\, c_1\exp\left(-c_2|x-y|\left(1+\log\frac{T|x-y|}{t}\right)^{(\b-1)/\b}\right) \text{ when } \,C_*T|x-y|\ge t.
\]
\end{prop}

\pf
We let $r:=|x-y|$ and fix $C_*\in(0,1)$. Note that $r\ge C_*^{-1}t/T>t/T\ge 1$ and $r\exp(-r^{\b})\le \exp(-1)(<1)$ for $\b> 1$. So  we only consider the case $Tr\exp(-r^{\b})<t$ $ (\le C_*Tr)$ which is equivalent to $r \left(\log(Tr/t)\right)^{-1/\beta}> 1$.
Let $k\ge 2$ be a positive integer such that
\begin{align}\label{e:k}
1<r \left(\log\frac{Tr}{t}\right)^{-1/\beta}\le k<r \left(\log\frac{Tr}{t}\right)^{-1/\beta}+1< 2r \left(\log\frac{Tr}{t}\right)^{-1/\beta}.
\end{align}
Then we have that
\begin{align}\label{e:t/k}
\frac{t}{k}\le \frac{t}{r}\left(\log\frac{Tr}{t}\right)^{1/\beta}\le T\cdot \sup_{s\ge C_*^{-1}}s^{-1}(\log s)^{1/\b}=:t_0<\infty
\end{align}

By our assumption on $D$, there is a length parameterized curve $l\subset D$ connecting $x$ and $y$ such that the total length $|l|$ of $l$ is less than or equal to $\lambda_1 r$ and $\delta_D(l(u))\ge \lambda_2 a\sqrt{t}$ for every $u\in [0, |l|]$.
We define
$r_t:=(2^{-1}{\lambda_2}a\sqrt{t})\wedge ((6\lambda_1)^{-1}(\log ({Tr}/{t}))^{1/\b})$.
Then by \eqref{e:k} and the assumption $\log (C_*^{-1})<\log(Tr/t)$, we have that
\begin{align}\label{e:r_t}
0<r_0:=\left(\frac{\lambda_2a\sqrt{T}}{2} \right)\wedge \left(\frac{(\log C_*^{-1})^{1/\b}}{6\lambda_1} \right) \le r_t\le\frac{1}{6\lambda_1}\left(\log\frac{Tr}{t} \right)^{1/\b}< \frac{r}{3\lambda_1 k}.
\end{align}

Define $x_i:=l(i |l|/k)$ and $B_i:=B(x_i, r_t)$ for $i=0,1,2, \ldots, k$ then  $\delta_D(x_i)\ge \lambda_2a\sqrt{t}>r_t$ and $B_i\subset D$. For every $y_i \in B_i$, we have that $\delta_D(y_i)\ge 2^{-1} \lambda_2 a \sqrt{t}>2^{-1} \lambda_2 a \sqrt{t/k}$ and
\begin{align}\label{e:y_i}
|y_i-y_{i+1}|\le|x_i-x_{i+1}|+2r_t \le\left(\lambda_1 +\frac{2}{3 \lambda_1} \right)\frac{r}{k}.
\end{align}
Thus by Proposition \ref{step1} and \ref{step2} along with the definition of $j$, \eqref{e:k}, \eqref{e:t/k} and  \eqref{e:y_i}, there exist constants $c_i>0, i=1,\ldots, 5$ such that
\begin{align}\label{p(t/k)}
&p_D(t/k, y_i, y_{i+1}) \ge c_{1}\left(\left(\frac{t}{k}\right)^{-d/2}\wedge\frac{t }{k}\cdot j(|y_i-y_{i+1}|)\right)\ge \, c_{2} \left(1 \wedge \Big(\frac{t}{k}\frac{e^{-c_{3}(r/k)^{\b}}}{\left(r/k\right)^{d+\alpha}}\Big)\right)\nn\\
& \ge \,c_{4} \frac{t}{Tr}\left(\frac{k}{r}\right)^{d+\alpha-1}e^{-c_{3}(r/k)^{\b}}\, \ge \, c_{4} \frac{t}{Tr}\left(\log \frac{Tr}{t}\right)^{-\frac{d+\alpha-1}{\b}}\left(\frac{t}{Tr}\right)^{c_{3}}\ge \, c_{4}\left(\frac{t}{Tr}\right)^{c_{5}}.
\end{align}

Therefore, by the semigroup property, \eqref{e:r_t} and \eqref{p(t/k)}, we conclude that
\begin{align*}
p_D(t, x, y)&\ge \int_{B_1}\cdots\int_{B_{k-1}} p_D(t/k, x, y_1)\cdots p_D(t/k, y_{k-1}, y ) dy_1\cdots dy_{k-1}\\
&\ge \left(c_4\left(\frac{t}{Tr}\right)^{c_5}\right)^k\Pi^{k-1}_{i=1}|B_i|\,\ge\, \left(\frac{c_6 t}{Tr}\right)^{c_5k}\\
&\ge c_7 \exp\left(-c_5k\left(\log\frac{Tr}{c_8t}\right)\right)\ge c_7 \exp\left(-c_9r\left(\log\frac{Tr}{t}\right)^{1-1/\b}\right)\\
&\ge c_7 \exp\left(-c_9r\left(1+\log\frac{Tr}{t}\right)^{1-1/\b}\right).
\end{align*}
\qed

\begin{prop}\label{step5}
Let $\b= \infty$ and $a, T$ and $C_*\in(1/2,1)$ be positive constants. Then there exist positive constants $c_i=c_i(a, T, C_*,\lambda_1, \lambda_2)$, $i=1,2$ such that  for every $t\in [T, \infty)$ and $x, y\in D$ with $\delta_D(x)\wedge \delta_D(y) \ge a\sqrt{ t}$, we have
\[
p_D(t,x,y) \,\ge \, c_1\exp\left(-c_2|x-y|\left(1+\log \frac{T|x-y|}{t}\right)\right) \text{ when }\, C_*T|x-y|\ge t.
\]
\end{prop}

\pf
Let $r:=|x-y|$ and fix $C_*\in (1/2,1)$. Since $T\le t \le C_*Tr$, we note that $1\le C_*r$.
By our assumption on $D$, there is a length parameterized curve $l\subset D$ connecting $x$ and $y$ such that the total length $|l|$ of $l$ is less than or equal to $\lambda_1 r$ and $\delta_D(l(u))\ge \lambda_2 a\sqrt{t}$ for every $u\in [0, |l|]$.
Let $k\ge 2 $ be a positive integer satisfying
\begin{align}\label{e:k2}
1<8 \lambda_1C_* r\leq k <8 \lambda_1C_* r+1\le (8 \lambda_1+1)C_* r.
\end{align}

Define $r_t:= (\lambda_2 a \sqrt{t}/2) \wedge 8^{-1}$, $x_i:=l(i|l|/k)$ and  $B_i:=B(x_i, r_t)$ for $i=0,1,\ldots,k$.
Then $\delta_D(x_i)> 2r_t$ and $B_i\subset B(x_i, 2 r_t)\subset D$.
For every $y_i\in B_i$, since $t/k<t/(8\lambda_1C_* r)\le T/(8\lambda_1)$, we have $\delta_D(y_i)>r_t>c_1 \sqrt{t/k}$ for some constant $c_1=c_1(a, T, \lambda_1, \lambda_2)>0$. Also, for each $y_i\in B_i$,
\begin{align}\label{e:y_i2}
|y_i-y_{i+1}|\leq |y_i-x_i|+|x_i-x_{i+1}|+|x_{i+1}-y_{i+1}| \leq
\frac{1}{8}+\frac{|l|}{k}+ \frac{1}{8}<\frac{\lambda_1 r}{8\lambda_1 C_*r}+\frac{1}{4}\le\frac{1}{2} .
\end{align}

By Proposition \ref{step1} and \ref{step2} along with the definition of $j$,  \eqref{e:y_i2} and the fact that $t/k<T/(8\lambda_1)$,  there are constants $c_i=c_i(a, T,  \lambda_1)>0$, $i=2,\ldots,4$, such that for $(y_i, y_{i+1})\in B_i\times B_{i+1}$,
\begin{align}\label{e:p(t/k)2}
p_D(t/k, y_i, y_{i+1})\ge c_2\left((t/k)^{-d/\alpha} \wedge \frac{t/k}{|y_i-y_{i+1}|^{d+\alpha}}\right)\ge  c_3 \left(1 \wedge t/k \right) \ge  c_4 \,t/(Tk).
\end{align}
Thus, by the semigroup property combining  the fact $r_t\ge r_T \wedge 8^{-1}$, \eqref{e:k2} and \eqref{e:p(t/k)2}, we obtain that
\begin{align*}
&p_D(t,x,y)\ge \int_{B_1}\ldots\int_{B_{k-1}}p_D(t/k,x,y_1)\ldots p_D(t/k, y_{k-1},y) dy_{k-1}\ldots dy_1\ge \left( \frac{c_4\,t}{Tk}\right)^k\Pi^{k-1}_{i=1} |B_i| \\
&\ge  \left(\frac{c_5\,t}{Tk}\right)^{k}\ge c_6 \left(\frac{c_7\,t}{Tr}\right)^{k} \ge c_6 \exp\left(-c_8r\log\frac{Tr}{c_7t}\right) \ge \exp\left(-c_9r\left(1+\log\frac{Tr}{t}\right)\right).
\end{align*}
\qed

Recall that  $B_R=B(x_0, R)$. Note that a exterior ball $\ol{B}^c_R$ is a domain in which the path distance is comparable to the Euclidean distance with characteristics $(\lambda_1, \lambda_2)$  independent of $x_0$ and $R$. Hence, the previous propositions yield the following Theorem.

\begin{thm}\label{T:int}
Let $a$ and $T$ be positive constants. Then for any $\b\in[0,\infty]$, there exists positive constants $c_i=c_i(a, \b, T)$ ($c=c(a)$ when $\b=0$, respectively) , $i=1,2$, such that  for every $R>0$, $t\in [T,\infty)$ ($t>0$ when $\beta=0$, respectively) and $x, y\in  \ol{B}_R^c$ with $\delta_{ \ol{B}_R^c}(x)\wedge \delta_{ \ol{B}_R^c}(y)\ge a\vp^{-1}(t)$, we have
$$
p_{\ol{B}_R^c}(t, x, y)\,\ge \, c_1
\Psi^2_{c_2, T}(t, |x-y|)
$$
where $\Psi^2_{c_2, T}(t, r)$ is defined in \eqref{eq:qd}.
\end{thm}
\pf
Let $r:=|x-y|$.  For any $\b\ge 0$, if $ \vp(r)<t$, by Proposition \ref{step1}, we have the conclusion.

Suppose $t\le\vp (r)$.
When $\b\in[0, 1]$, we have the conclusion by Proposition \ref{step2} and Proposition \ref{step3}.
When $\b\in (1, \infty)$, using Proposition \ref{step3} and Proposition \ref{step4} , and when $\b=\infty$, using Proposition \ref{step3} and Proposition \ref{step5}, we have the conclusion.
\qed

\section{Lower bound estimates}

In this section, we assume that the dimension $d>2\cdot{\bf 1}_{\{ \beta \in (0, \infty]\}}   + \alpha\cdot {\bf 1}_{\{ \beta =0\}}$.
To establish the lower bound estimates in Theorem \ref{T:main}(2)--(4), we first consider the lower bound estimates on $p_{\ol{B}_R^c}(t, x, y)$ for  $t\ge T$ ($t>0$ when $\b=0$, respectively) where $B_R$ is a ball of radius $R>0$ centered at  $x_0$. Since all following estimates are independent of $x_0$, we may assume that $x_0=0$.

We define the Green function $G(x, y)$ of $Y$ in $\R^d$ as $G(x, y):=\int_0^{\infty} p(t, x, y)dt$  for every $x, y\in \R^d$.
Then by the fact that
$\int_0^{\infty}(t^{-d/\alpha}\wedge tr^{-d-\alpha})dt \asymp r^{\alpha-d}$ for $d>\alpha$ when $\b=0$ and by \cite [Theorem 6.1] {CKK3} when $\b\in(0, \infty]$,
 we have that
\begin{align}\label{e:green}
G(x, y)\asymp  \left(|x-y|^{\alpha-d}+|x-y|^{2-d}\cdot {\bf 1}_{\{\b\in(0, \infty]\}}\right).
\end{align}

For any Borel set $A\subset \R^d$, define the first exit time of $A$ as $\tau_A=\inf\{t>0: Y_t\notin A\}$   and the first hitting time of $A$ as $T_A=\inf\{t>0: Y_t\in A\}$. The next lemma provide us the beginning point for the lower bound estimates which proof is almost identical to that of \cite[ Lemma 4.1]{CKS7} using \eqref{e:green}, so we omit the proof.

\begin{lemma}\label{L:4.1}
There is a constant $C_3>1$ such that for all $R>0$,
\begin{align*}
\qquad C_3^{-1}\frac{R^{d}}{R^{\alpha}+R^2}&\left(|x|^{\alpha-d}+|x|^{2-d}\cdot {\bf 1}_{\{\b\in(0, \infty]\}}\right)\le \,\P_x(T_{\ol{B}_R}<\infty)\nn\\
&\qquad\qquad\le C_3 \,\frac{R^{d}}{R^{\alpha}+R^2}\left(|x|^{\alpha-d}+|x|^{2-d}\cdot {\bf 1}_{\{\b\in(0, \infty]\}}\right), \qquad  \text{ for }\,\, |x|\ge 2R.
\end{align*}
\end{lemma}

\medskip

The following ideas of obtaining the lower bound estimates on $p_{\ol{B}_R^c}(t, x, y)$  are motivated by that of  Section 5 in  \cite{CKS7} and for the sake of completeness, we give proofs detail.
For the simplicity of the notation, hereafter  for any $ y\in \R^d\backslash \{0\}$ and $r>0$, we define
$H(y,r):= \{z\in B(y, r):z\cdot y\ge0\}.$
Recall that  $\vp(r)=r^2\cdot {\bf 1}_{\{\beta\in (0, \infty]\}}+r^{\alpha}\cdot{\bf 1}_{\{\beta=0\}}$ and $\vp^{-1}(t)=t^{1/2}\cdot {\bf 1}_{\{\beta\in (0, \infty]\}}+t^{1/\alpha}\cdot{\bf 1}_{\{\beta=0\}}$.

\begin{lemma}\label{L:5.1}
Let $T$ be a positive constant. Then for any  $\b\in [0, \infty]$,  there exists  constants $\varepsilon=\varepsilon(\b, T)>0$ and  $M_1=M_1(\b,  T)\ge 3$ ($\varepsilon>0$ and $M_1\ge 3$ when $\b=0$, respectively) such that the following holds: for any $R>0$, $t\in [T, \infty)$  ($t>0$ when $\b=0$, respectively) and $x, y$ satisfying $|x|>M_1R$, $|y|>R$ and $y\in B(x, 9\vp^{-1}(t))$, we have
$$
\P_x\left( Y_t^{\ol {B}_R^c}\in H(y, \vp^{-1}(t)/2)\right)\ge \varepsilon. $$
\end{lemma}

\pf
Applying  \eqref{eq:laT} (Applying \eqref{eq:smT} and  \eqref{eq:laT} when $\b=0$, respectively) and by the change of variable with $v=z/\vp^{-1}({t})$, for any $t\ge T$ ($t>0$ when $\b=0$, respectively),   there are constants $c_i=c_i(\b, T)>0$ ($c_i>0$ when $\b=0$, respectively), $i=1,\cdots,3$ such that
\begin{align*}
\P_x&\left(Y_t\in H(y, \vp^{-1}({t})/2)\right)\ge  \inf_{w\in B(y, 9\vp^{-1}(t))}\P_w\left(Y_t\in H(y, \vp^{-1}(t)/2)\right)\\
&\ge\, c_1 \inf_{w\in B(y, 9\vp^{-1}(t))} \int_{H(y, \vp^{-1}(t)/2)} \Psi^2_{C_1, T}(t,|w-z|)\cdot{\bf 1}_{\{\beta\in (0, \infty]\}}+\left(t^{-d/\alpha}\wedge t|w-z|^{-d-\alpha}\right)\cdot {\bf 1}_{\{\beta=0\}}\, dz\\
&\ge\, c_ 2\inf_{w\in B(y, 9\vp^{-1}(t))} \int_{H(y, \vp^{-1}(t)/2)} \frac{1}{\vp^{-d}(t)}\left(\exp \left(-C_1\frac{|w-z|^2}{t}\right)\cdot {\bf 1}_{\{\beta\in (0, \infty]\}}+{\bf 1}_{\{\beta=0\}}\right) \,dz\\
&= \, c_3 \inf_{w_0\in B(y_0, 9)} \int_{H(y_0, 1/2)} \exp \left(-C_1 |w_0-v|^2\right)\cdot {\bf 1}_{\{\beta\in (0, \infty]\}}+{\bf 1}_{\{\beta=0\}}\, dv\\
&\ge \,2^{-1}c_3 |B(0, 1/2)|\left(e^{-C_1 10^2}\cdot {\bf 1}_{\{\beta\in (0, \infty]\}}+{\bf 1}_{\{\beta=0\}}\right)
\end{align*}
where $y_0:=y/\vp^{-1}(t)$ and $w_0:=w/\vp^{-1}(t)$. When $\b=0$,  since $|w-z|\le 10t^{1/\alpha}$, the third inequality holds.
Hence, there is $\varepsilon \in (0, 1/4)$ so that for any $t\ge T$ ($t>0$ when $\b=0$, respectively), $x\in \R^d$ and $y\in B(x, 9\vp^{-1}(t))$, we have
\begin{align}\label{e:P1}
\varepsilon<\frac{1}{2}\P_x\left(Y_t\in H(y, \vp^{-1}(t)/2)\right).
\end{align}

For $d>2\cdot{\bf 1}_{\{ \beta \in (0, \infty]\}}+ \alpha\cdot {\bf 1}_{\{ \beta =0\}}$ and the constant $C_3>1$  in Lemma \ref{L:4.1}, we may choose $M_1\ge 3$ so that $C_3(M_1^{2-d}+M_1^{\alpha-d}\cdot{\bf1}_{\{\b\in(0, \infty]\}})\le \varepsilon$. For any $x$ with $|x|>M_1R$, by Lemma \ref{L:4.1}, we have that
\begin{align}\label{e:P2}
\P_x\left(\tau_{\ol{B}^c_R}\le t\right)= \P_x\left(T_{\ol{B}_R}<\infty\right)&\le C_3\frac{R^d}{R^{\alpha}+R^2}(|x|^{2-d}+|x|^{\alpha-d}\cdot{\bf1}_{\{\b\in(0, \infty]\}})\nn\\
&\le C_3  \left(\frac{R^2}{R^{\alpha}+R^2}M_1^{2-d}+\frac{R^{\alpha}}{R^{\alpha}+R^2}M_1^{\alpha-d}\cdot{\bf1}_{\{\b\in(0, \infty]\}}\right)\nn\\
&\le C_3(M_1^{2-d}+M_1^{\alpha-d}\cdot{\bf1}_{\{\b\in(0, \infty]\}})\le \varepsilon.
\end{align}
Hence, combining \eqref{e:P1} and \eqref{e:P2}, we obtain that
\begin{align*}
\P_x\left(Y_t^{\ol{B}^c_R}\in H(y, \vp^{-1}(t)/2)\right)= &\,\P_x\left(\tau_{\ol{B}^c_R}>t\right) -\P_x\left(Y_t^{\ol{B}^c_R}\notin H(y, \vp^{-1}(t)/2); \tau_{\ol{B}^c_R}>t\right)\\
\ge&\, \P_x\left(\tau_{\ol{B}^c_R}>t\right)-\P_x\left(Y_t\notin H(y, \vp^{-1}(t)/2)\right)\\
\ge&\, (1-\varepsilon)-(1-2\varepsilon)=\varepsilon.
\end{align*}
\qed

\begin{lemma}\label{L:5.2}
Let $T>0$, $\b\in [0, \infty]$, and $M_1=M_1( \b, T/8)\ge 3$ ($ M_1\ge 3$ when $\b=0$, respectively) be the constant in Lemma \ref{L:5.1}. Then there exists a positive constant $c=c(\b, T)>0$ ($c>0$ when $\b=0$, respectively) such that for any $R>0$, $t\in[T, \infty)$ ($t>0$ when $\b=0$, respectively) and $x, y$ satisfying $|x|>M_1R$, $|y|>M_1 R$ and  $|x-y|\le \vp^{-1}(t)/6$, we have that
$p_{\ol{B}_R^c}(t, x, y)\ge c/\vp^{-d}(t).$
\end{lemma}

\pf
Without loss of generality we may assume that $|y|\ge |x|$. If $\delta_{\ol{B}_R^c}(y)>\vp^{-1}(t)/2$, then $\delta_{\ol{B}_R^c}(x)\ge \delta_{\ol{B}_R^c}(y)-|x-y|\ge \vp^{-1}(t)/3$, and hence the lemma follows immediately from Proposition \ref{step1}.

Now we assume that $\delta_{\ol{B}_R^c}(y)\le \vp^{-1}(t)$/2. By the semigroup property and the parabolic Harnack inequality(see \cite [Theorem 4.11]{CKK3}) , we have
\begin{align}\label{e:1}
p_{\ol{B}_R^c}(t, x, y)&\ge \int_{H(y, \vp^{-1}(t/2))} p_{\ol{B}_R^c}(t/2,x,z) p_{\ol{B}_R^c}(t/2, z, y)dz\nn\\
&\ge c_1 \P_x\left(Y_{t/2}^{\ol{B}_R^c}\in H(y,\vp^{-1}(t/2))\right) p_{\ol{B}_R^c}\left(t/2-\vp(2\delta_{\ol{B}_R^c}(y))/4, y, y\right).
\end{align}

Note that $t\ge s:= t/2-\vp(2\delta_{\ol{B}_R^c}(y))/4\ge t/4\ge T/4$ ($s\ge t/4>0$ when $\b=0$, respectively). So by the semigroup property, the Cauchy-Schwarz inequality and Lemma \ref{L:5.1}, we obtain that
\begin{align}\label{e:2}
p_{\ol{B}_R^c}(s, y, y)\ge& \int_{H(y, \vp^{-1}(s)/2)} \left(p_{\ol{B}_R^c}(s/2, y,z)\right)^2 dz\nn\\
\ge & \frac{2}{|B(y, \vp^{-1}(s)/2)|}\P_y\left(Y^{\ol{B}_R^c}_{s/2}\in H(y, \vp^{-1}(s)/2)\right)^2\ge c_2 /\vp^{-d}(s)\ge c_2  /\vp^{-d}(t).
\end{align}
Applying Lemma \ref{L:5.1} again and \eqref{e:2} to \eqref{e:1}, we have that
$p_{\ol{B}_R^c}(t, x, y)\ge c_3 /\vp^{-d}(t)$.
\qed

\begin{prop}\label{P:5.4}
Let $T>0$, $\b\in [0, \infty]$, and $M_1=M_1(\b,T/16)\ge 3$ ($M_1\ge 3$ when $\b=0$, respectively) be the constant in Lemma \ref{L:5.1}. Then there exist positive constants $c=c(\b, T)$ and $C_4=C_4(\b, T)$ ($c, \,C_4>0$ when $\b=0$, respectively) such that for any $R>0$, $t\in[T, \infty)$ ($t>0$ when $\b=0$, respectively) and $x, y$ satisfying $|x|>M_1R$, $|y|>M_1 R$,  we have that
$
p_{\ol{B}_R^c}(t, x, y)\ge c \Psi^2_{C_4, T}(t, |x-y|),
$
where $\Psi^2_{a, T}(t, r)$ is defined in \eqref{eq:qd}.
\end{prop}

\pf
By Lemma \ref{L:5.2}, we only need to prove the proposition for $|x-y|>\vp^{-1}(t)/6$.

If $t/2\le \vp(60R)$, then $\delta_{\ol{B}_R^c}(x)\wedge \delta_{\ol{B}_R^c}(y)\ge (M_1-1)R\ge 2R\ge(30)^{-1}\vp^{-1}(t/2)$. In this case the Proposition holds by Theorem \ref{T:int}. So  we only consider the following case: $t\ge  T\wedge 2\vp(60R)$ ($t\ge 2\vp(60R)$ when $\b=0$, respectively) and $|x-y|>\vp^{-1}(t)/6$ . Without loss of generality, we may assume that $|y|\ge |x-y|/2$. Let $x_1:=x+20^{-1}\vp^{-1}(t/2)x/|x|$ then we have $B(x_1, 20^{-1}\vp^{-1}(t/2))\subset \ol{B}_{|x|}^c\subset \ol{B}_R^c$.

For every $z\in B(x_1, 20^{-1}\vp^{-1}(t/2))$, we obtain
\begin{align}\label{inB}
|x-z|\le \frac {1}{20}\vp^{-1}(t/2)+|x_1-z|\le \frac{1}{10}\vp^{-1}(t/2)\le \frac{1}{6}\vp^{-1}(t/2).
\end{align}
Since $|y|\ge |x-y|/2$ and $R\le 60^{-1}\vp^{-1}(t/2)$, we have
\begin{align}\label{del_1}
\delta_{\ol{B}_R^c}(y)=|y|-R\ge\frac{1}{2}|x-y|-\frac{1}{60}\vp^{-1}(t/2)> \frac{1}{12}\vp^{-1}(t)-\frac{1}{60}\vp^{-1}(t/2)\ge \frac{1}{15}\vp^{-1}(t/2).
\end{align}
For $z\in B(x_1, 60^{-1}\vp^{-1}(t/2))$, we have
\begin{align}\label{del_2}
\delta_{\ol{B}_R^c}(z)=|z|-R&\ge |x_1|-|x_1-z|-\frac{1}{60}\vp^{-1}(t/2)\nn\\
&\ge |x|+\frac{1}{20}\vp^{-1}(t/2)-\frac{1}{60}\vp^{-1}(t/2)-\frac{1}{60}\vp^{-1}(t/2)\ge \frac{1}{60}\vp^{-1}(t/2)
\end{align}
and
$$|z-y|\le |z-x|+|x-y|\le \frac{1}{15} \vp^{-1}(t/2)+|x-y|\le 2|x-y|.$$
By the semigroup property, Lemma \ref{L:5.2} with \eqref{inB},  Theorem \ref{T:int} with \eqref{del_1} and \eqref{del_2} and the fact $r\to \Psi^2_{a, T}(t,r)$ is decreasing,  there exist constants $c_i=c_i(\b, T)>0$ ($c_i>0$ when $\b=0$, respectively), $i=1, \ldots, 4$ such that
\begin{align*}
p_{\ol{B}_R^c}(t, x, y)=&\int_{\ol{B}_R^c} p_{\ol{B}_R^c}(t/2, x, z)p_{\ol{B}_R^c}(t/2, z, y)dz\\
\ge& \int_{B(x_1, \vp^{-1}(t/2)/60)}p_{\ol{B}_R^c}(t/2, x, z)p_{\ol{B}_R^c}(t/2, z, y)dz\\
\ge &\, c_1\, \int_ {B(x_1, \vp^{-1}(t/2)/60)} 1/(\vp^{-d}(t/2))\Psi^2_{c_2, T/2}(t/2, |z-y|)dz\\
\ge &\,c_3 \,\Psi^2_{2c_2, T}(t, 2|x-y|)\ge \,c_4 \,\Psi^2_{2^3c_2 , T}(t, |x-y|).
\end{align*}
The last inequality holds by (1) in Lemma \ref{L:p2} with $a=2c_2$ and $b=2$ and we have proved the proposition.
\qed

The following elementary lemma is used to prove the lower bound estimates on $p_D(t, x, y)$ where $t\in [T, \infty)$ ($t>0$ when $\b=0$, respectively). Recall the function $\Psi^1_{a, b, T}(t, r)$   which is defined in \eqref{eq:qd}.

\begin{lemma}\label{L:com}
Let $ K, R , b$ and $t_0$ be fixed positive constants and $\b\in[0,\infty]$. Suppose that $x, x_1\in \R^d$ satisfy $|x-x_1|=K^2R$. Then there exists a positive constant $ c=c( K, R, b, t_0, \b)$ such that for any $a>0$ and $z\in \R^d$, we have
$\Psi^1_{a, b, t_0}\left(t_0, 5|x-z|/4\right)\ge c~ \Psi^1_{a, b, t_0}\left(t_0, 2|x_1-z|\right).$
\end{lemma}

\pf
Let $r:=|x-z|$ and $r_1:=|x_1-z|$. For any $z\in B(x, KR)\cup B(x_1, KR)$, we have that $r\le (K+1)KR$. So  $\Psi^1_{a, b, t_0}( t_0, 5r/4)$ is bounded below and the lemma holds.

Suppose that  $z\notin B(x, KR)\cup B(x_1, KR)$.  When $r\le 4K^2R\vee 4/5$, then $\Psi^1_{a, b,  t_0}( t_0, 5r/4)$ is bounded below and hence the lemma holds. Let   $r> 4K^2R\vee 4/5$. By the triangle inequality, we have that  $3r/4<r-K^2R\le r_1\le r+K^2R<5r/4$ and hence $1\le 5r/4\le 5r_1/3\le 2r_1$. In this case, since $r\to \Psi^1_{a, b, t_0}(t_0, r) $ is non-increasing, the lemma holds.
\qed

Now, we are ready to prove the lower bound estimates on $p_D(t, x, y)$.
For the remainder of this paper, we assume that  $\eta\in(\alpha/2,1]$ and  $D$ is an exterior $C^{1,\eta}$ open set in $\R^d$ with $C^{1, \eta}$ characteristics $(r_0, \Lambda_0)$ and $D^c\subset B(0, R)$ for some $R>0$. Such an open set $D$ can be disconnected. When $\b\in(1, \infty]$ and  $|x-y|\ge 4/5$, we will consider the following two cases that $x, y$ are in the same component and in different components in $D$, separately.

\medskip

\noindent {\bf Proof of Theorem \ref{T:main}(2)--(3)}
 Due to Theorem \ref{T:sm}(4) and the domain monotonicity of $p_D(t, x, y)$, the Theorem holds when $x, y$ are in the same bounded connected component of $D$.
So we only need to prove Theorem \ref{T:main}(2)--(3.a).

When $\b=0$, by Theorem \ref{T:sm}(1), we may assume that $t\ge T$. Without loss of generality, we may assume that $T=3$. For $x$ and $y$ in $D$, let $v\in \R^d$ be any unit vector satisfying $x\cdot v\ge 0$ and $y\cdot v\ge 0$. Let $M_2:=M_1(\b, 3(16)^{-1})(\ge 3)$, where $M_1$ is the constant in  Lemma \ref{L:5.1}. Define
$$x_1:=x+M_2^2 R v\,\, \text{ and } \,\, y_1:=y+M_2^2Rv.$$

By the semigroup property and Theorem \ref{T:sm}(1)-(2), we have that for every $t-2\ge 1$ and $x, y\in D$,
\begin{align}\label{e:lp2}
p_D(t, x, y)&=\int_D\int_Dp_D(1, x,z)p_D(t-2, z, w)p_D(1, w, y)dzdw\nn\\
&\ge c_1 (1\wedge \delta_D(x))^{\alpha/2} (1\wedge \delta_D(y))^{\alpha/2}f_2(t, x, y),
\end{align}
where  $C_2$ and $\gamma$ are given constants in Theorem \ref{T:sm}  and
\begin{align}\label{e:f_2}
f_2(t, x, y)=&\int_{B(0, M_2R)^c\times B(0, M_2R)^c} (1\wedge \delta_D(z))^{\alpha/2}\Psi^1_{C_2, \gamma, 1}(1, 5|x-z|/4)\nn\\
&\qquad\cdot \,p_D(t-2, z, w) (1\wedge \delta_D(w))^{\alpha/2}\Psi^1_{C_2, \gamma, 1}(1, 5|y-w|/4) dzdw.
\end{align}

Let $A_2:=\max\{(C_1C_4)^{1/{(\b\wedge 1)}}, 2\gamma^{2/\b}, 2C_1C_2\}(\ge 2)$ ($A_2= 2$ when $\b=0$, respectively) where  $C_1$ is the constant in \eqref{eq:smT}, \eqref{eq:laT} and $C_4$ is the constant in Proposition \ref{P:5.4}. By  Lemma \ref{L:com}  and \eqref{eq:smT}, there exists $c_i=c_i(\b)>0, i=2, \ldots, 4$ such that
\begin{align}\label{e:5r/4}
\Psi^1_{C_2, \gamma, 1}(1, 5|x-z|/4)&\ge c_2\Psi^1_{C_2, \gamma, 1}(1, 2|x_1-z|) \nn\\
&\ge c_3 \Psi^1_{C_1^{-1}, \gamma^{-1}, 1}(1, A_2|x_1-z|)\ge c_4 p(1, A_2x_1, A_2z)~\mbox{ and }\nn\\
\Psi^1_{C_2, \gamma, 1}(1, 5|y-w|/4)&\ge c_2\Psi^1_{C_2, \gamma, 1}(1, 2|y_1-w|) \nn\\
&\ge c_3 \Psi^1_{C_1^{-1}, \gamma^{-1}, 1}(1, A_2|y_1-w|)\ge c_4 p(1, A_2y_1, A_2w).
\end{align}
When $\b\in(0, \, \infty]$, the second inequalities hold by (2) in Lemma \ref{L:p1}
 along with $t_0=1$, $a=C_1$, $b=C_2$, $c=\gamma$, $N_1=2$ and  $N_2=A_2$ and the fact $A_2\ge  2( C_1C_2\vee \gamma^{2/\b})$. When $\b=0$, the second inequalities hold since $A_2= 2$.

For $z, w\in B(0, M_2R)^{c}$ and $t-2\in [1, \infty)$, by Proposition \ref{P:5.4} and \eqref{eq:laT}, we have that
\begin{align}\label{e:t-2}
p_D(t-2, z, w)&\ge p_{\ol{B}_R^c}(t-2, z, w)\ge c_5\Psi^2_{C_4, 1}(t-2, |z-w|)\nn\\
&\ge c_6 \Psi^2_{C_1^{-1}, 1}(t-2, A_2|z-w|)\ge c_7 p(t-2, A_2z, A_2w).
\end{align}
For the third inequality above, we use (2) in Lemma \ref{L:p2} along with $T=1$, $a=C_4$, $b=C_1$ and $N=A_2$ and the fact $A_2\ge (C_1C_4)^{1/(\b\wedge 1)}$ when $\b\in(0, \infty]$. When $\b=0$, the third inequality holds since $A_2\ge 1$.

For $z\in B(0, M_2R)^c$, $\delta_D(z)\ge\delta_{\ol{B}^c_{R}}(z)= |z|-R\ge M_2R-R$.
So  applying  \eqref{e:5r/4} and \eqref{e:t-2} to \eqref{e:f_2} and by the change of variables $\hat{z}=A_2z$, $\hat{w}=A_2w$ and semigroup property,  we have that
\begin{align}\label{in:f_2}
f_2(t, x, y)&\ge c_{8} \int_{B(0, M_2R)^c\times B(0, M_2R)^c} p(1, A_2x_1, A_2z)p(t-2, A_2z, A_2w)p(1, A_2y_1, A_2w) dzdw\nn\\
&\ge c_{9} \int_{B(0, A_2M_2R)^c\times B(0, A_2M_2R)^c} p_{B(0, A_2M_2R)^c}(1, A_2x_1, \hat{z}) p_{B(0, A_2M_2R)^c}(t-2, \hat{z}, \hat{w})\nn\\
&\qquad\qquad\qquad\qquad\cdot p_{B(0, A_2M_2R)^c}(1, A_2y_1, \hat{w})d\hat{z}d\hat{w}\nn\\
&=c_{9}\, p_{B(0, A_2M_2R)^c}(t, A_2x_1, A_2y_1).
\end{align}
Since $A_2|x_1|\wedge A_2|y_1|\ge M_2(A_2M_2R)$, by Proposition \ref{P:5.4} and (1) in Lemma \ref{L:p2} with $a=C_4$ and $b=A_2$, we have that
\begin{align}\label{in:r}
p_{B(0, A_2M_2R)^c}(t, A_2x_1, A_2y_1)\ge& c_{10} \,\Psi^2_{C_4, T}(t, A_2|x_1-y_1|)\nn\\
=& c_{10} \Psi^2_{C_4, T}(t, A_2|x-y|)\ge c_{11}\Psi^2_{A_2^2C_4, T}(t, |x-y|).
\end{align}
Combining \eqref{e:lp2} with \eqref{in:f_2} and \eqref{in:r}, we have proved the lower bound estimates in Theorem \ref{T:main}(2)--(3.a).
\qed

For the remainder of this section, we assume that $T>0$, $\b\in(1, \infty)$ and $(t, x, y)\in [T, \infty)\times D\times D$ where $|x-y|\ge 4/5$ and $x, y$ are in different components of $D$.

It is clear that there exists $0<\kappa \le 1/2$ which is depending on $\Lambda_0$ and $d$ such that  for all $x\in \overline{D}$ and  $r\in (0,r_0]$ there is a ball $B(A_r(x), \kappa r)\subset D\cap B(x,r)$.
Hereinafter,  we assume  that $A_r(x)$ is such the point in $D$.

\begin{lemma}\label{bD}
Suppose that $D_b\subset B(0, R)$ be a bounded connected component of $D$.
Then there exists a positive constant $c=c(\b, \eta,  r_0, \Lambda_0, T)$ such that  for every $t\ge T$ and $x\in D_b$, we can find a ball $B\subset D_b$ such that
\begin{align*}
\int_{B} p_{D_b}(2^{-1}t-3^{-1}T, x, z)dz\ge c\,e^{-t \lambda^{D_b}}  \delta_{D_b}(x)^{\alpha/2}
\end{align*}
where $-\lambda^{D_b}<0$ be the largest eigenvalue of the generator of $Y^{D_b}$.
\end{lemma}

\pf
For any $x\in D_b$, let $z_x\in \ol{D_b}$ be the point so that $|z_x-x|=\delta_{D_b}(x)$.
Let $x_1:=A_{r_0}(z_x)$ and $B:=B(x_1, \kappa r_0)$. For any $z\in B$, we have that $\delta_{D_b}(z)\ge \kappa r_0$.
Hence since $2^{-1}t-3^{-1}T\ge6^{-1}T$, by Theorem \ref{T:sm}(4) along with the bounded connected component $D_b$,  there exist constants $c_i=c_i(\b, \eta, r_0, \Lambda_0, T)>0$, $i=1,\ldots,3$ such that for any $x\in D_b$
\begin{align*}
\int_{B}p_{D_b}(2^{-1}t-3^{-1}T, x, z)dz&\ge c_1e^{-t\lambda^{D_b}}\int_{B}\delta_{D_b}(x)^{\alpha/2}\delta_{D_b}(z)^{\alpha/2}dz\nn\\
&\ge c_2 e^{-t\lambda^{D_b}}\delta_{D_b}(x)^{\alpha/2}\int_{B}dz\ge c_3e^{-t\lambda^{D_b}}\delta_{D_b}(x)^{\alpha/2}.
\end{align*}
\qed

Now, we are ready to prove the lower bound estimates on $p_D(t, x, y)$ for any  $\b\in(1, \infty)$ and $(t, x, y)\in [T, \infty)\times D\times D$ where $|x-y|\ge 4/5$ and $x,y$ are in different components of $D$.

\bigskip

\noindent {\bf Proof of  Theorem \ref{T:main}(4)}
Let $D(x)$ and $D(y)$ be connected components containing $x$ and $y$, respectively with $D(x)\cap D(y)\ne\emptyset$. Without loss of generality, we may assume that $D(x)$ is a bounded connected component and  $T=3$.

By the semigroup property and  the domain monotonicity of $p_D(t, x, y)$, we first observe that
\begin{align}\label{eq:dc}
p_D(t, x, y)\ge \int_{D(x)}\int_{D(y)}p_{D(x)}(2^{-1}t-1, x, z)p_{D}(2, z, w)p_{D(y)}(2^{-1}t-1, y, w)dwdz.
\end{align}

For bounded connected component $D_j$ of $D$ and  the largest eigenvalue $-\lambda_{j}<0$ of the generator $Y^{D_j}$, define  $\ol{\lambda}:=\max\{\lambda_j:j=1,\ldots,n\}$.
By Lemma \ref{bD},  there exist a ball $B_x\subset D(x)$ and a constant $c_1=c_1(\b, \eta, r_0, \Lambda_0)>0$ such that
\begin{align}\label{D(x)}
\int_{B_x}p_{D(x)}(2^{-1}t-1, x, z)dz\ge c_1e^{-t\ol{\lambda}} \delta_{D}(x)^{\alpha/2}.
\end{align}

Similarly, if $D(y)$ is a bounded connected component,  we have that
$\int_{B_y}p_{D(y)}(2^{-1}t-1, y, w)dw\ge c_2e^{-t\ol{\lambda}} \delta_{D}(y)^{\alpha/2}$ for some  a ball $B_y\subset D(y)$ and a constant $c_2>0$.
For any $(z, w)\in  B_{x}\times B_{y}$, note that $r_0\le |z-w|\le 2R$ and $\delta_{D}(z)\wedge \delta_{D}(w)\ge c_3$. So  by Theorem \ref{T:sm}(1) and (3),  we have that
$ \inf_{(z, w)\in B_{x}\times B_{y}}p_{D}(2, z, w)\ge c_4$.
Hence, we have the conclusion when $D(x)$ and $D(y)$ are bounded connected components of $D$.

When $D(y)$ is an unbounded connected component, let $y_1:=y+2Ry/|y|$ and $B_{y_1}:=B(y_1, 2^{-1}R ) $ $\subset D(y)$. For any $w\in B_{y_1}$, we have that $\delta_{D(y)}(w)\ge R/2$ and $|y-w|\le |y-y_1|+|y_1-w|\le 5R/2$.
Hence for $2^{-1}t-1\ge 1/2$, by Theorem \ref{T:main}(2)--(3.a) and the fact $ t/2-1\asymp t$, there exist constants $c_i=c_i(\b, \eta, r_0, \Lambda_0 , R)>0$, $i=5,\ldots,8$ such that
\begin{align}\label{D(y)}
\int_{B_{y_1}}p_{D(y)}(2^{-1}t-1, y, w)dw&\ge c_5\int_{B_{y_1}}(1\wedge \delta_{D(y)}(y))^{\alpha/2}(1\wedge \delta_{D(y)}(w))^{\alpha/2}t^{-d/2}\exp(-c_6|y-w|^2/t)dw\nn\\
&\ge c_7 (1\wedge \delta_{D}(y))^{\alpha/2}t^{-d/2}\int_{B_{y_1}}dw=c_8 (1\wedge \delta_{D}(y))^{\alpha/2}t^{-d/2}.
\end{align}
For any $(z, w)\in B_{x}\times B_{y_1}$, we have that $\delta_{D}(z)\wedge \delta_{D}(w)\ge c_9$ and
$$|z-w|\le |z-x|+|x-y|+|y-w|\le 2R+|x-y|+5R/2\le c_{10}|x-y|.$$
The last inequality holds since $|x-y|\ge 4/5$.
So  by Theorem \ref{T:sm}(1) and (3), there are constants $c_i=c_i(\b, \eta, r_0, \Lambda_0 , R)>0$, $i=11,\ldots,14$  such that
\begin{align}\label{infP}
\inf_{(z, w)\in B_{x}\times B_{y_1}}p_D(2, z, w)\ge c_{11}\left(\frac{\exp(-c_{12}|z-w|^{\b})}{|z-w|^{d+\alpha}}\wedge 1\right)\ge c_{13}\frac{\exp(-c_{14}|x-y|^{\b})}{|x-y|^{d+\alpha}}.
\end{align}
Combining \eqref{D(x)}, \eqref{D(y)} and \eqref{infP} with \eqref{eq:dc}, we have the conclusion when $D(x)$  is a bounded connected component and $D(y)$ is an unbounded connected component of $D$.
\qed

\begin{remark}
Let $D$  be an exterior $C^{1, \eta}$ open set in $\R^d$ with $C^{1, \eta}$ characteristics $(r_0, \Lambda_0)$ and $D^{\,c}\subset B(0, R)$. Then the number of bounded connected components of $D$ is finite, say $D_1,\ldots, D_n$.
According to the proof of Theorem \ref{T:main}(4), there exists a constant $c>0$ such that  if $x, y\in D$ are in different bounded connected components of $D$
\begin{align*}
p_D(t, x, y)\ge c \,\delta_D(x)^{\alpha/2}\delta_D(y)^{\alpha/2}\exp\left(-t\sum_{j=1}^n \lambda_j\Big({\bf 1}_{D_j}(x)+{\bf 1}_{D_j}(y)\Big)\right)
\end{align*}
where $-\lambda_j<0 $ is the largest eigenvalue of the generator $Y^{D_j}$, $j=1,\cdots,n$.
\end{remark}

\section{Green function estimate}
In this section, we present a proof of Corollary \ref{C:green}. We recall that $G_D(x,y)=\int_0^\infty p_D(t, x, y)dt$.
When $\b=0$, the proof of Corollary \ref{C:green} is similar to that of   \cite[Corollary 1.5]{CT}, we only consider the case $\b\in(0, \infty]$.
\medskip

\noindent {\bf Proof of Corollary \ref{C:green}}
By Corollary \ref{Coro}, there exist constants $c_i>1$, $i=1,2$ such that
\begin{align*}
G_D(x, y)&\le c_1 \int_0^1 \left(1\wedge \frac{\delta_D(x)}{t^{1/\alpha}}\right)^{\alpha/2} \left(1\wedge \frac{\delta_D(y)}{t^{1/\alpha}}\right)^{\alpha/2}\Psi^1_{c_2^{-1}, \gamma^{-1}, 30}(t, |x-y|/6)dt\nn\\
&\qquad+ c_1 \left(1\wedge \delta_D(x)\right)^{\alpha/2}\left(1\wedge \delta_D(y)\right)^{\alpha/2}\int_1^{\infty}\Psi^2_{c_2^{-1}, 1}(t, |x-y|) dt \qquad\qquad\mbox{ and }\nn\\
G_D(x, y)&\ge c_1^{-1}\cdot{\bf 1}_{\{|x-y|<4/5\}} \int_0^1 \left(1\wedge \frac{\delta_D(x)}{t^{1/\alpha}}\right)^{\alpha/2} \left(1\wedge \frac{\delta_D(y)}{t^{1/\alpha}}\right)^{\alpha/2}\Psi^1_{c_2, \gamma, 30}(t, 5|x-y|/4)dt\nn\\
&\qquad+ c_1^{-1}\cdot{\bf 1}_{\{|x-y|\ge4/5\}} \left(1\wedge \delta_D(x)\right)^{\alpha/2}\left(1\wedge \delta_D(y)\right)^{\alpha/2}\int_1^{\infty}\Psi^2_{c_2, 1}(t, |x-y|)dt
\end{align*}
where $\gamma$ is the constant in Theorem \ref{T:sm}.

Without loss of generality, we may assume that $c_2=1$ and we define $I_1$, $I_2$ and $II$ by
\begin{align*}
&I_1:=\int_0^1 \left(1\wedge \frac{\delta_D(x)}{t^{1/\alpha}}\right)^{\alpha/2} \left(1\wedge \frac{\delta_D(y)}{t^{1/\alpha}}\right)^{\alpha/2}\left(t^{-d/\alpha}\wedge t|x-y|^{-\alpha-d}\right) dt\nn\\
&I_2:=\int_0^1 \left(1\wedge \frac{\delta_D(x)}{t^{1/\alpha}}\right)^{\alpha/2} \left(1\wedge \frac{\delta_D(y)}{t^{1/\alpha}}\right)^{\alpha/2}\Psi^1_{1, \gamma^{-1}, 30}(t, |x-y|/6)dt\nn\\
&II:= \int_1^\infty \Psi^2_{1, 1}(t, |x-y|)dt.
\end{align*}
For any $a, b>0$, if $b|x-y|<1$, we have that $\Psi^1_{1, a, 30}(t, b|x-y|)\asymp t^{-d/\alpha}\wedge t|x-y|^{-\alpha-d}$. So  when $|x-y|< 4/5$, it suffices to show that
\begin{align}\label{eq:rl}
&I_1\asymp \left(\frac{1}{|x-y|^{d-\alpha}}+\frac{1}{|x-y|^{d-2}}\right)\left(1\wedge \frac{\delta_D(x)}{|x-y|}\right)^{\alpha/2}\left(1\wedge \frac{\delta_D(y)}{|x-y|}\right)^{\alpha/2}\mbox{ and }\nn\\
&II\le c_3\le c_4\left(\frac{1}{|x-y|^{d-\alpha}}+\frac{1}{|x-y|^{d-2}}\right).
\end{align}
When $|x-y|\ge 4/5$, we will show that
\begin{align}
&I_2\le c_5\left(\frac{1}{|x-y|^{d-\alpha}}+\frac{1}{|x-y|^{d-2}}\right)\left(1\wedge \delta_D(x)\right)^{\alpha/2}\left(1\wedge \delta_D(y)\right)^{\alpha/2}\mbox{ and}\label{eq:ru}\\
&II\asymp\frac{1}{|x-y|^{d-2}}\asymp\left(\frac{1}{|x-y|^{d-\alpha}}+\frac{1}{|x-y|^{d-2}}\right).\label{eq:ru1}
\end{align}

Let $r:=|x-y|$.
Suppose that  $r<  4/5$. By [\cite{CKS},  (4.3), (4.4) and (4.6)], we have 
\begin{align}\label{I_1}
I_1&\asymp \frac{1}{r^{d-\alpha}}\left(1\wedge \frac{\delta_D(x)}{r}\right)^{\alpha/2}\left(1\wedge \frac{\delta_D(y)}{r}\right)^{\alpha/2}\nn\\
&\asymp \left(\frac{1}{r^{d-\alpha}}+\frac{1}{r^{d-2}}\right)\left(1\wedge \frac{\delta_D(x)}{r}\right)^{\alpha/2}\left(1\wedge \frac{\delta_D(y)}{r}\right)^{\alpha/2}.
\end{align}
Note that for every $ s \in[0, \infty]$,
\begin{align}\label{e:int}
\int_s^{\infty} t^{-d/2}e^{-r^2/t}dt =r^{2-d}\int_0^{r^2/s} u^{d/2-2}e^{-u}du.
\end{align}
For $r<1$ and $1< t$, we have  $\Psi^2_{1,1}(t, r)=t^{-d/2}e^{-r^2/t}$ and 
\begin{align}\label{II}
II=r^{2-d}\int_0^{r^2} u^{d/2-2}e^{-u}du\asymp r^{2-d}\int_0^{r^2} u^{d/2-2}du=\frac{2}{d-2}.
\end{align}
Hence we obtain \eqref{eq:rl} by \eqref{I_1} and \eqref{II}.

Suppose that $r\ge 4/5$. Note that for $0< t\le1$,  we have
\begin{align*}
\Psi^1_{1,\gamma^{-1}, 30}(t, r/6)=
&\begin{cases}
t^{-d/\alpha}\wedge t(r/6)^{-d-\alpha}e^{-\gamma^{-1}(r/6)^{\b}}&\le t(r/6)^{-d-\alpha}\qquad\mbox{ for } \b\in(0,1]\nn\\
t\exp(-((r/6)(\log(5r/t))^{(\b-1)/\b}\wedge (r/6)^\b))&\le te^{-c_6r}\qquad\qquad\mbox{ for } \b\in(1, \infty)\nn\\
\left(t/(5r)\right)^{r/6}&\le t^{2/15}e^{-c_6r}\,\,\qquad\mbox{ for } \b=\infty
\end{cases}\nn\\
\le& ~c_7t^{2/15}r^{-d-\alpha}
\end{align*}
for some constant $c_i=c_i(\b)>0$, $i=6, 7$. Thus by the change of variable $u=r^\alpha/t$, there exist constants $c_i>0, i=8,9$ such that
\begin{align}\label{I_2}
I_2\le &  ~c_7r^{-d-\alpha}\int_0^1 t^{\frac{2}{ 15}} \left(1\wedge \frac{\delta_D(x)}{t^{1/\alpha}}\right)^{\alpha/2}\left(1\wedge \frac{\delta_D(y)}{t^{1/\alpha}}\right)^{\alpha/2}dt\nn\\
= &~ c_7r^{-d+\frac{2}{ 15}\alpha} \int_{r^{\alpha}}^\infty u^{-\frac{2}{ 15}-2} \left(1\wedge \frac{\sqrt{u}\delta_D(x)^{\alpha/2}}{r^{\alpha/2}}\right)\left(1\wedge \frac{\sqrt{u}\delta_D(y)^{\alpha/2}}{r^{\alpha/2}}\right) du\nn\\
= &~ c_7r^{-d+\frac{2}{ 15}\alpha} \int_{r^{\alpha}}^\infty u^{-\frac{2}{ 15}-1} \left(\frac{1}{\sqrt{u}}\wedge \frac{\delta_D(x)^{\alpha/2}}{r^{\alpha/2}}\right)\left(\frac{1}{\sqrt{u}}\wedge \frac{\delta_D(y)^{\alpha/2}}{r^{\alpha/2}}\right) du\nn\\
\le &~c_8r^{-d+\frac{2}{ 15}\alpha} \int_{r^{\alpha}}^\infty u^{-\frac{2}{ 15}-1} du \left(1 \wedge \frac{\delta_D(x)}{r}\right)^{\alpha/2}\left(1\wedge \frac{\delta_D(y)}{r}\right)^{\alpha/2} \nn\\
=&~\frac{15}{2}c_8r^{-d} \left(1 \wedge \frac{\delta_D(x)}{r}\right)^{\alpha/2}\left(1\wedge \frac{\delta_D(y)}{r}\right)^{\alpha/2}
\le~c_9r^{2-d} \left(1 \wedge \delta_D(x)\right)^{\alpha/2}\left(1\wedge \delta_D(y)\right)^{\alpha/2}
\end{align}
and it yields \eqref{eq:ru}. For \eqref{eq:ru1},  because of \eqref{II}, we may assume that  $r\ge1$. By \eqref{e:int}, we have that
\begin{align*}
II&\ge \int_1^\infty t^{-d/2}e^{-r^2/t} dt\ge r^{2-d}\int_0^1 u^{d/2-2}e^{-u}du\ge c_{10} r^{2-d}\qquad\mbox{ and }\\
II&\le\int_1^{r^{2-(\b\wedge1)}} t^{-d/2}e^{-r^{(\b\wedge 1)}}dt +\int_{r^{2-(\b\wedge1)}}^\infty t^{-d/2}e^{-r^2/t} dt\\
&\le c_{11} e^{-r^{(\b\wedge 1)}}+r^{2-d}\int_0^{r^{(\b\wedge 1)}} u^{d/2-2}e^{-u}du\le c_{12} r^{2-d}.
\end{align*}
This implies $ II\asymp r^{2-d}$ and hence \eqref{eq:ru1} holds.
So we have proved the Corollary.
\qed

\bigskip

\noindent {\bf Acknowledgement}
This paper is a part of the author's Ph.D. thesis. She thanks Professor Panki Kim, her Ph.D. thesis advisor, for his guidance and encouragement.

\begin{small}

\end{small}

\vskip 0.3truein

{\bf Kyung-Youn Kim}

Department of Mathematical Sciences,
Seoul National University,
Building 27, 1 Gwanak-ro, Gwanak-gu,
Seoul 151-747, Republic of Korea

E-mail: \texttt{eunicekim@snu.ac.kr}


\begin{thebibliography}{99}

\bibitem{BBCK}   {M.~T. Barlow,  R.~F. Bass,  Z.-Q. Chen and M. Kassmann,}
Non-local Dirichlet forms and symmetric jump processes.
\textit{Trans. Amer. Math. Soc.} \textbf{361} (2009), 1963--1999.

\bibitem{BGR1} {K. Bogdan, T. Grzywny and M. Ryznar,}
Heat kernel estimates for the fractional Laplacian with Dirichlet conditions.
\textit{Ann. Probab.} \textbf{38} (2010), 1901-–1923.

\bibitem{BGR2} {K. Bogdan, T. Grzywny and M. Ryznar,}
Density and tails of unimodal convolution semigroups.
\textit{J. Funct. Anal.} \textbf{266(6)} (2014), 3543 – 3571.


\bibitem{CKK}  Z.-Q. Chen, P. Kim and T. Kumagai,
Weighted Poincar\'e inequality and heat kernel estimates for finite range jump processes.
\textit{Math. Ann.} \textbf{342(4)} (2008), 833--883.

\bibitem{CKK2}  Z.-Q. Chen, P. Kim and  T. Kumagai,
On heat kernel estimates and parabolic Harnack inequality for jump processes on metric measure spaces.
\textit{ Acta Math. Sin. (Engl. Ser.)} \textbf{25} (2009), 1067--1086.

\bibitem{CKK3}  Z.-Q. Chen, P. Kim and  T. Kumagai,
Global heat kernel estimates for symmetric jump processes.
\textit{Trans. Amer. Math. Soc.} \textbf{363(9)} (2011), 5021--5055.

\bibitem{CKS}  Z.-Q. Chen, P. Kim, and R. Song,
Heat kernel estimates for Dirichlet fractional Laplacian.
\textit{J. European Math. Soc.} \textbf{12} (2010), 1307--1329.

\bibitem{CKS2} Z.-Q. Chen, P. Kim and R. Song,
Dirichlet heat kernel estimates for $\Delta^{\alpha/2} + \Delta^{\beta/2}$.
\textit{Ill. J. Math.} \textbf{54} (2010), 1357--1392.

\bibitem{CKS3} Z.-Q. Chen, P. Kim and R. Song,
Sharp heat kernel estimates for relativistic stable processes in open sets.
\textit{Ann. Probab.} \textbf{40} (2012), 213--244.

\bibitem{CKS6}  Z.-Q. Chen, P. Kim and R. Song,
Global heat kernel estimates for relativistic stable processes in half-space-like open sets.
\textit{Potential Anal.} \textbf{36} (2012), 235--261.

\bibitem{CKS7}  Z.-Q. Chen, P. Kim and R. Song,
Global heat kernel estimates for relativistic stable processes in exterior open sets.
\textit{J. Funct. Anal.} \textbf{263(2)} (2012), 448-475.

\bibitem{CKS8} Z.-Q. Chen, P. Kim and R. Song,
Dirichlet heat kernel estimates for subordinate Brownian motions with Gaussian components.
\textit{J. Reine Angew. Math.} (2014) in press.


\bibitem{CKS9} Z.-Q. Chen, P. Kim and R. Song,
Dirichlet heat kernel estimates for rotationally symmetric L\'evy processes.
\textit{Proc. Lond. Math. Soc. (3)}  \textbf{109} (2014), 90--120.


\bibitem{CK}  Z.-Q. Chen and  T. Kumagai,
Heat kernel estimates for stable-like processes on d-sets.
\textit{Stochastic Process Appl.} \textbf{108} (2003), 27--62.

\bibitem{CK2}  Z.-Q. Chen and  T. Kumagai,
Heat kernel estimates for jump processes of mixed types on metric measure spaces.
\textit{Probab. Theory Relat. Fields} \textbf{140} (2008), 277--317.

\bibitem{CT} Z.-Q. Chen and J. Tokle,
Global heat kernel estimates for fractional Laplacian in unbounded open sets.
\textit{Probab. Theory Relat. Fields} \textbf{149} (2011), 373--395.


\bibitem{CZ} K.~L. Chung and Z. Zhao,
{\it From Brownian Motion to Schr\"{o}dinger's Equation}.
Springer, Berlin, 1995.

\bibitem{FOT} M. Fukushima, Y. Oshima and M. Takeda,
{\it Dirichlet Forms and Symmetric Markov Processes}.
Walter De Gruyter, Berlin, 1994.

\bibitem{KaSz1} K. Kaleta and P. Sztonyk,
Upper estimates of transition densities for stable-dominated semigroups.
\textit{J. Evol. Equ.} \textbf{13(3)} (2013), 633--650.

\bibitem{KaSz2} K. Kaleta and P. Sztonyk,
Estimates of transition densities and their derivatives for jump L\'evy processes.
2013,  arXiv:1307.1302v2.


\bibitem{KK} K. Kim and P. Kim,
Two-sided estimates for the transition densities of symmetric Markov processes dominated by stable-like processes in $C^{1,\eta}$ open sets. 
\textit{Stochastic Process Appl.} \textbf{124(9)} (2014), 3055--3083.



\bibitem{Sz1} P. Sztonyk,
Approximation of stable-dominated semigroups.
\textit{Potential Anal.} \textbf{33} (2010), 211--226.


\end{thebibliography}
\end{document}